\documentclass[12pt]{preprint}

\usepackage{amsmath,amssymb}
\usepackage{mhequ,mhenvs}
\usepackage{graphicx}
\usepackage{hyperref}
\usepackage{booktabs}
\usepackage{todo}
\usepackage{times}
\newtheorem{assumption}[lemma]{Assumption}

\newtheorem{defn}[lemma]{Definition}

\renewcommand{\Omega}{D}
\newcommand{\pn}{P^{\mbox \tiny{N}}}

\newcommand{\muy}{\mu}
\newcommand{\llambda}{\theta}

\newcommand{\zi}{z^{(i)}}

\newcommand{\vi}{v^{(i)}}
\newcommand{\zo}{z^{(1)}}
\newcommand{\uo}{u_1}
\newcommand{\vo}{v^{(1)}}
\newcommand{\zt}{z^{(2)}}
\newcommand{\ut}{u_2}
\newcommand{\vt}{v^{(2)}}
\newcommand{\bbE}{\mathbb{E}}

\newcommand{\bbP}{\mathbb{P}}
\newcommand{\bbR}{\mathbb{R}}

\newcommand{\bbK}{\mathbb{K}}
\newcommand{\bbZ}{\mathbb{Z}}

\newcommand{\dhh}{d_{\mbox {\tiny{\rm Hell}}}}
\newcommand{\dtv}{d_{\mbox {\tiny{\rm TV}}}}

\newcommand{\igmas}{{\cal C}}

\newcommand{\ud}{\mathrm{d}}
\newcommand{\PN}{\mathbb{P}^{N}}
\newcommand{\vN}{v^{N}}
\newcommand{\lap}{\Delta}
\newcommand{\grad}{\nabla}
\newcommand{\cdiv}{\nabla\cdot}
\newcommand{\ddt}{\frac{\ud}{\ud t}}
\newcommand{\lt}{\left}

\newcommand{\nx}{\|_X}

\newcommand{\ny}{\|_Y}

\newcommand{\G}{\mathcal{G}}

\newcommand{\cC}{\mathcal{C}}
\newcommand{\cH}{\mathcal{H}}

\newcommand{\kk}{\tiny{\mathrm K}}
\newcommand{\jj}{\tiny{\rm J}}
\newcommand{\kb}{{\mbox K}}

\newcommand{\jb}{{\mbox J}}
\newcommand{\cHs}{\mathcal{H}^s}

\newcommand{\eps}{{\epsilon}}

\newcommand{\bpf}{\begin{proof}}
\newcommand{\epf}{\end{proof}}

\newcommand{\cA}{\mathcal{A}}

\newcommand{\cG}{\mathcal{G}}

\newcommand{\cN}{\mathcal{N}}

\newcommand\scal[2][ ]{\ifthenelse{\equal{#1}{ }}{\langle#2\rangle}{}
        \ifthenelse{\equal{#1}{b}}{\bigl\langle#2\bigr\rangle}{}
        \ifthenelse{\equal{#1}{B}}{\Bigl\langle#2\Bigr\rangle}{}
        \ifthenelse{\equal{#1}{bb}}{\biggl\langle#2\biggr\rangle}{}
        \ifthenelse{\equal{#1}{BB}}{\Biggl\langle#2\Biggr\rangle}{}}

\let\E\expect

\def\L^#1{\mathrm{L}^{\!#1}}

\newcommand{\bbT}{\mathbb{T}}

\renewcommand{\epsilon}{\varepsilon}

\renewcommand{\le}{\leqslant}
\renewcommand{\ge}{\geqslant}

\renewcommand{\theenumi}{\roman{enumi}}
\renewcommand{\theenumiii}{\arabic{enumiii}}

\newcommand{\be}{\begin{eqnarray}}
\newcommand{\ee}{\end{eqnarray}}
\renewcommand{\E}{\bbE}

\begin{document}

\title{Approximation of Bayesian\\ Inverse Problems for PDEs}
\author{S.L. Cotter, M. Dashti and A.M. Stuart \\
          Mathematics Institute \\
University of Warwick \\
Coventry CV4 7AL, UK}
\maketitle
\begin{abstract}
Inverse problems are often ill-posed,                                 
with solutions that depend sensitively  
on data. In any numerical approach to the
solution of such problems, regularization 
of some form is                   
needed to counteract the resulting instability. This paper 
is based on an approach to regularization, 
employing a Bayesian formulation of the 
problem, which leads to a notion of                                   
well-posedness for inverse problems, 
at the level of probability measures. 

The stability which results from this                                 
well-posedness may be used as the 
basis for quantifying the approximation,
in finite dimensional spaces,
of inverse problems for functions. 
This paper contains a theory which 
utilizes the stability to estimate
the distance between the true and
approximate posterior distributions,
in the Hellinger metric, in terms of
error estimates for approximation of
the underlying forward problem. 
This is potentially useful as it allows
for the transfer of estimates from
the numerical analysis of forward problems
into estimates for the solution of the
related inverse problem.
In particular controlling
differences in the Hellinger metric 
leads to 
control on the differences between expected 
values of polynomially bounded functions and
operators, including the mean and covariance
operator.
                                                                      
The ideas are illustrated with the  
classical inverse problem for the   
heat equation, and then applied to some
more complicated non-Gaussian inverse problems 
arising in data assimilation, involving determination of
the initial condition for the Stokes
or Navier-Stokes equation from Lagrangian
and Eulerian observations respectively.
\end{abstract} 

\section{Introduction}

In applications it is frequently of interest to solve 
{\em inverse problems} \cite{kir96,Tar05}:
to find $u$, an input to a mathematical model, given
$y$ an observation of (some components of, or functions of) 
the solution of the model. We have an equation of the form
\begin{equation}
\label{eq:basic}
y=\cG(u)
\end{equation}
to solve for $u \in X$, given $y \in Y$, where $X,Y$ are
Banach spaces.
We refer to evaluating $\cG$ as solving the
{\em forward problem}\footnote{In the applications in
this paper $\cG$ is found from composition of the
forward model with some form of observation operator,
such as pointwise evaluation at a finite set of
points. The resulting observation operator
is often denoted with the letter $\cH$ in the atmospheric sciences
community \cite{ICGL97}; 
because we need $\cH$ for Hilbert space later on, 
we use the symbol $\cG$.}. We refer to $y$ as {\em data} or
{\em observations}.
It is typical of inverse problems that 
they are {\em ill-posed}\index{ill-posed}:
there may be no solution, or the solution may not be unique and may depend
sensitively on $y$. 
For this reason some form of regularization is often employed
\cite{ehn96} to stabilize computational approximations.

We adopt
a Bayesian approach to regularization \cite{bs94} which 
leads to the notion of
finding a {\em probability measure} $\muy$ on $X$, containing information
about the relative probability of different states $u$, given the
data $y$. 
Adopting this approach is natural in situations where 
an analysis of the 
source of data reveals that the observations $y$
are subject to noise. A more appropriate model 
equation is then often of the form
\begin{equation}
\label{eq:basicn}
y=\cG(u)+\eta
\end{equation}
where $\eta$ is a mean-zero random variable, whose statistical properties
we might know, or make a reasonable mathematical model
for, but whose actual value is unknown to us; we refer to $\eta$
as the {\em observational noise}.
We assume that it is possible to
describe our prior knowledge
about $u$, before acquiring data,
in terms of a {\em prior probability measure} 
$\mu_0$. It is then possible to use Bayes' formula to calculate the 
{\em posterior probability measure} $\muy$ for
$u$ given $y$.

In the infinite dimensional setting the most
natural version of Bayes theorem is a statement
that the posterior measure is absolutely
continuous with respect to the prior \cite{Stuart10}
and that the Radon-Nikodym derivative (density) between them
is determined by the data likelihood. 
This gives rise to the formula
\begin{equation}
\label{eq:cs}
\frac{d\muy}{d\mu_0}(u)=\frac{1}{Z(y)}\exp\bigl(-\Phi(u;y)\bigr)
\end{equation}
where the {\em normalization constant}
$Z(y)$ is chosen so that $\muy$
is a probability measure:
\begin{equation}
\label{eq:Rnormz}
Z(y)=\int_{X} \exp\bigl(-\Phi(u;y)\bigr) d\mu_0(u).
\end{equation}
In the case
where $y$ is finite dimensional and $\eta$ has
Lebesegue density $\rho$ this is simply 
\begin{equation}
\label{eq:irnd}
\frac{d\muy}{d\mu_0}(u) \propto \rho(y-\cG(u)).
\end{equation}
More generally $\Phi$ is determined by the distribution of
$y$ given $u.$
We call $\Phi(u;y)$ the {\em potential},\index{potential}
and sometimes, for brevity, refer to evaluation of $\Phi(u;y)$ 
for a particular $u \in X$,
as {\em solving the forward problem} as it is defined
through $\cG(\cdot).$ 
Note that the solution to the inverse problem is
a probability measure $\muy$ which is defined
through a combination of solution of the forward
problem $\cG$, the data $y$ and a prior probability measure $\mu_0.$

In general it is hard to obtain information from a formula
such as \eqref{eq:cs} for a probability
measure. One useful approach to extracting
information is to use {\em sampling}:
generate a set of points
$\{u^{(k)}\}_{k=1}^K$ distributed (perhaps only approximately)
according to $\muy.$ In this context it is noteworthy
that the integral $Z(y)$ appearing in formula \eqref{eq:cs}
is not needed to enable implementation of MCMC
methods to sample from the desired measure.  
These methods incur an error which is well
understood and which decays as $\sqrt{K}$ \cite{Liu01}. 
However for inverse problems on function space
there is a second source of error, arising from the need
to approximate the inverse problem in a finite dimensional
subspace of dimension $N$. 
The purpose of this paper is to quantify such
approximation errors. The {\em key idea}
is that we transfer approximation
properties of the forward problem $\Phi$ into approximation
properties of the inverse problem defined by \eqref{eq:cs}.

Since the solution to the Bayesian inverse problem
is a probability measure we will need to use metrics
on probability measures to quantify the effect of
approximation. We will employ the Hellinger metric
$\dhh$ from Definition \ref{def:tvd} because this 
leads directly to bounds on the 
approximation error incurred when calculating
the expectation of functions. This property is summarized in
Lemma \ref{l:tvhell2}. Combining these ideas we will
find that finite dimensional approximation 
leads to
an error in the calculation of expectation of functions
which tends to zero as $\psi(N)$ tends to infinity, for some function $\psi(N)$
{\em determined by approximation of the forward problem.} 

In section \ref{sec:frame} we provide the general
approximation theory, for measures $\mu$ given by
\eqref{eq:cs}, upon which the remainder of the
paper builds. Section \ref{sec:heat} employs
this approximation theory to study the classical
inverse problem of determining the initial condition
for the heat equation from observation of the
solution at a positive time. In section \ref{sec:lag}
we study the inverse problem of determining the initial condition
for the Stokes equation, given a finite set of
observations of Lagrangian
trajectories defined through the time-dependent velocity
field solving the Stokes equation; this section
also includes numerical results showing the convergence
of the posterior distribution under refinement
of the finite dimensional approximation, as predicted
by the theory.  Section \ref{sec:eul}
is devoted to the related inverse problem 
of determining the initial condition
for the Navier-Stokes equation, given direct observation
of the time-dependent velocity field at a finite set of
points at positive times.

A classical approach to the regularization of
inverse problems is through the least squares approach and
Tikhonov regularization \cite{ehn96,Tar05}; a good
overview of this approach, in the context of data assimilation
problems in fluid mechanics such as those studied
in sections \ref{sec:lag} and \ref{sec:eul}, is \cite{Nic02}
and the connection between the least squares and
Bayesian approaches for applications in fluid mechanics
is overviewed in \cite{ApteJV08}.
The Bayesian formulation to inverse problems in general
is overviewed in the text \cite{ks04}. Note,
however, that the methodology employed there
is typically one in which the problem is first
discretized, and then ideas from Bayesian
statistics are applied to the resulting
finite dimensional problem. The approach taken
in this paper is to first formulate the Bayesian
inverse problem {\em on function space} and then
study approximation. As
in many areas of applied mathematics -- for example,
optimal control -- formulation of the problem in function
space, followed by discretization will lead to better
algorithms and better understanding. This approach 
is laid out conceptually in \cite{Tar05} for inverse
problems, but the underlying mathematics is not
developed, except for some particular linear and Gaussian
problems. Indeed, for linear problems, the
Bayesian approach on function space may be found 
in an early paper of Franklin \cite{Fr70}, including
study of the heat equation, the subject of section \ref{sec:heat}.
More recently there has been some work on
finite dimensional linear inverse problems, using
the Bayesian approach to regularization,
and considering infinite dimensional limits \cite{hopi, NP08}
and in the limit of disappearing observational noise
\cite{hopi2}. A general approach to the formulation,
and well-posedness, of inverse problems, 
adopting a Bayesian approach on function
space, is undertaken in \cite{CDRS08}; furthermore
applications to problems in fluid mechanics are
given in that paper and we will build on
this material in sections \ref{sec:lag} and \ref{sec:eul}.

\section{General Framework}
\label{sec:frame}

In this section we establish three useful results which concern
the effect of approximation on the posterior probability
measure $\muy$ given by \eqref{eq:cs}. 
These three results are Theorem \ref{t:wp2}, Corollary 
\ref{c:wp2} and Theorem \ref{t:wp2f}. The key point to
notice about these results is that they simply require the
proof of various bounds and approximation properties
{\em for the forward problem}, and yet they yield approximation
results concerning the Bayesian inverse problem. The connection
to probability comes only through the choice of the space
$X$, in which the bounds and approximation properties must
be proved, which must have full measure under the prior $\mu_0.$

The probability measure of interest 
\eqref{eq:cs} is defined through a 
density with respect to a prior 
reference measure $\mu_0$ which, by shift of origin,
we take to have mean zero. Furthermore, 
we assume that this reference measure is Gaussian
with covariance operator $\cC.$ 
We write $\mu_0=\cN(0,\cC).$ In fact we only use the Fernique
Theorem \ref{thm:fern} for $\mu_0$ and the results
may be trivially extended to all measures which satisfy
the conclusion of this theorem. The Fernique Theorem 
holds for all Gaussian measures on a separable Banach
space \cite{bog98}, and
also for other measures with tails which decay at least
as fast as a Gaussian.

It is demonstrated in \cite{Stuart10}
that in many applications, including those
considered here, the potential $\Phi(\cdot;y)$ 
satisfies certain natural bounds
on a Banach space
$\Bigl(X,\|\cdot\nx\Bigr)$, 
contained in the original Hilbert space on which
$\mu_0$ is defined, and of full measure under $\mu_0$
so that $\mu_0(X)=1.$ Such bounds are summarized in the
following assumptions. We assume that the data $y$ lies
in a Banach space $\Bigl(Y,\|\cdot\ny\Bigr).$
The key point about the form of Assumption \ref{ass:1}(i) is that
it allows use of the Fernique Theorem to control integrals
against $\mu.$ The assumption (ii) may be used to obtain
lower bounds on the normalization constant $Z(y).$ 

\begin{assumption}
\label{ass:1}
For some Banach space
$X$ with $\mu_0(X)=1,$
the function $\Phi:X\times Y \to \bbR$ 
satisfies the following: 
\begin{enumerate}
\item for every $\epsilon>0$ and $r>0$
there is $M=M(\epsilon,r) \in \bbR$ 
such that, for all $u \in X$ and $y\in Y$ with $\|y\ny<r$,
$$\Phi(u;y) \ge M-\epsilon\|u\nx^2;$$
\item for every $r>0$ there is a $L=L(r)>0$
such that, for all $u \in X$ and $y\in Y$ with 
$\max\{\|u\nx,\|y\ny\}<r$,
$$\Phi(u;y) \le L(r).$$
\end{enumerate}
\end{assumption}

For Bayesian inverse problems in which a finite
number of observations are made and the observation
error $\eta$ is mean zero Gaussian, the 
potential\index{potential} $\Phi$ has the form
\begin{equation}
\Phi(u;y)=\frac12|y-\cG(u)|_{\Gamma}^2
\label{eq:potf2}
\end{equation}
where $y\in\bbR^m$ is the data, $\cG:X \to \bbR^m$
is the forward model 
and $|\cdot|_{\Gamma}$ is
a covariance weighted norm on $\bbR^m$ given by
$|\cdot|_{\Gamma}=|\Gamma^{-\frac12}\cdot|$ and
$|\cdot|$ denotes the standard Euclidean norm.
In this case it is natural to express conditions
on the measure $\muy$ in terms of $\cG.$

\begin{assumption}
\label{ass:g1}
For some Banach space $X$ with $\mu_0(X)=1$,
the function $\cG:X \to \bbR^m$ 
satisfies the following: for
every $\epsilon>0$ there is $M=M(\epsilon)\in \bbR$ 
such that, for all $u \in X,$
$$|\cG(u)|_{\Gamma} \le \exp\bigl(\eps\|u\nx^2+M).$$
\end{assumption}

\begin{lemma}
\label{l:qwp}
Assume that $\Phi:X\times\bbR^m \to \bbR$ is given by
\eqref{eq:potf2} and let $\cG$ satisfy
Assumptions \ref{ass:g1}. 
Assume also that $\mu_0$ is a Gaussian measure satisfying
$\mu_0(X)=1.$
Then $\Phi$ satisfies Assumptions \ref{ass:1}. 
\end{lemma}

\begin{proof} 
Assumption \ref{ass:1}(i) is automatic 
since $\Phi$ is positive; assumption (ii)
follows from the bound 
$$\Phi(u;y) \le |y|^2_{\Gamma}+|\cG(u)|_{\Gamma}^2$$
and use of the exponential bound on $\cG$.
\end{proof}

Since the dependence on $y$ is not relevant 
we suppress it notationally and study measures $\mu$ given by
\begin{equation}
\label{eq:csagain2}
\frac{d\mu}{d\mu_0}(u)=\frac{1}{Z}\exp\bigl(-\Phi(u)\bigr)
\end{equation}
where the normalization constant\index{normalization
constant} $Z$ is given by 
\begin{equation}
\label{eq:snormz2}
Z=\int_{X} \exp\bigl(-\Phi(u)\bigr) d\mu_0(u).
\end{equation}
We approximate $\mu$ by approximating $\Phi$ over some
$N-$dimensional subspace of $X$. In
particular we define $\mu^N$ by
\begin{equation}
\label{eq:csagain23}
\frac{d\mu^N}{d\mu_0}(u)=\frac{1}{Z^N}\exp\bigl(-\Phi^N(u)\bigr)
\end{equation}
where 
\begin{equation}
\label{eq:Anormz2}
Z^N=\int_{X} \exp\bigl(-\Phi^N(u)\bigr) d\mu_0(u).
\end{equation}
The potential $\Phi^N$ should be viewed as
resulting from an approximation to the solution of the forward
problem.
Our interest is in translating approximation results
for $\Phi$ into approximation results for $\mu.$

The following theorem proves such a result, bounding 
the Hellinger distance, and hence by \eqref{eq:tvhrel}
the total variation distance,
between measures $\mu$ and $\mu^N$, in terms of the error
in approximating $\Phi.$ Again the particular exponential
dependence of the error constant for the forward
approximation is required so that we may use the Fernique
Theorem to control certain expectations arising in the
analysis.

\begin{theorem}
\label{t:wp2}
Assume that $\Phi$ and $\Phi^N$ satisfy 
Assumptions \ref{ass:1}(i),(ii) with constants uniform in $N$.
Assume also that, for any $\epsilon>0,$ there is $K=K(\epsilon)>0$
such that
\begin{equation}
\label{eq:forward}
|\Phi(u)-\Phi^N(u)|\le K\exp\bigl(\epsilon\|u\nx^2\bigr)\psi(N)
\end{equation}
where $\psi(N) \to 0$ as $N \to \infty$. Then
the measures $\mu$ and $\mu^N$ are close 
with respect to the Hellinger distance: there is a constant
$C$, independent of $N$, and such that
\begin{equation}
\label{eq:error}
\dhh(\mu,\mu^{N}) \le C\psi(N).
\end{equation}
Consequently all moments of $\|u\nx$ are 
${\cal O}\bigl(\psi(N)\bigr)$ close.
In particular the mean and, in the case $X$ is a Hilbert
space, the covariance operator, are 
${\cal O}\bigl(\psi(N)\bigr)$ close.
\end{theorem}

\begin{proof}
Throughout the proof, all integrals are over $X$.
The constant $C$ may depend upon $r$ and changes
from occurrence to occurrence.
Using Assumption \ref{ass:1}(ii) gives
$$|Z| \ge \int_{\{\|u\nx \le r\}} \exp\bigl(-L(r)\bigr)d\mu_0(u)
\ge \exp\bigl(-L(r)\bigr)\mu_0\{\|u\nx \le r\}.$$
This lower bound is positive because
$\mu_0$ has full measure on $X$ and is Gaussian so that all
balls in $X$ have positive probability.
We have an analogous lower bound for $|Z^N|.$

From Assumptions \ref{ass:1}(i) and \eqref{eq:forward}, 
using the fact that $\mu_0$ is a Gaussian probability measure 
so that the Fernique Theorem \ref{thm:fern} applies,
\begin{align*}
|Z-Z^N| &\le \int K\psi(N)\exp\bigl(\epsilon\|u\nx^2-M\bigr)
\exp\bigl(\epsilon\|u\nx^2\bigr)d\mu_0(u)\\
& \le C\psi(N).
\end{align*}

From the definition of Hellinger distance we have
\begin{align*}
2\dhh(\mu,\mu^N)^2 & = \int\Bigl(
Z^{-\frac12}\exp\bigl(-\frac{1}{2}\Phi(u)\bigr)-(Z^N)^{-\frac12}
\exp\bigl(-\frac{1}{2} \Phi^N(u)\bigr)\Bigr)^2d\mu_0(u)\\
&\le I_1+I_2 
\end{align*}
where
\begin{align*}
I_1&=\frac{2}{Z}\int\Bigl(
\exp\bigl(-\frac{1}{2}\Phi(u)\bigr)-\exp(-\frac{1}{2}\Phi^N(u)\bigr)\Bigr)^2
d\mu_0(u),
\\
I_2&=2|Z^{-\frac12}-(Z^N)^{-\frac12}|^2\int \exp(-\Phi^N(u)\bigr)
d\mu_0(u).
\end{align*}

Now, again using Assumptions \ref{ass:1}(i) and 
equation \eqref{eq:forward}, together with the
Fernique Theorem \ref{thm:fern}, 
\begin{align*}
\frac{Z}{2}I_1 & \le 
 \int \frac14 K^2\psi(N)^2\exp\bigl(3\epsilon\|u\nx^2-M\bigr)d\mu_0(u)\\
&\le C\psi(N)^2. 
\end{align*}

Note that the bounds on $Z,Z^N$ from below are
independent of $N$. Furthermore,
$$\int \exp\bigl(-\Phi^N(u)\bigr)d\mu_0(u) \le \int \exp\bigl(
\epsilon\|u\nx^2-M\bigr)d\mu_0(u)$$
with bound independent of $N$, by the Fernique 
Theorem \ref{thm:fern}. Thus
\begin{align*}
I_2 & \le C\bigl(Z^{-3} \vee (Z^N)^{-3}\bigr)|Z-Z^N|^2\\
&\le C\psi(N)^2. 
\end{align*}
Combining gives the desired continuity result in the
Hellinger metric.

Finally all moments of $u$ in $X$ are finite under  
the Gaussian measure $\mu_0$ by the Fernique 
Theorem \ref{thm:fern}. It follows
that all moments are finite under $\mu$ and $\mu^N$ 
because, for $f:X \to Z$ polynomially bounded, 
$$\bbE^{\mu} \|f\| \le \bigl(\bbE^{\mu_0} \|f\|^2\bigr)^{\frac12}
\bigl(\bbE^{\mu_0} \exp(-2\Phi(u;y))\bigr)^{\frac12}$$
and the first term on the right hand side is finite  since
all moments are finite under $\mu_0$,
whilst the second term may be seen to be finite by use of
Assumption \ref{ass:1}(i) and the Fernique Theorem \ref{thm:fern}.
\end{proof}

For Bayesian inverse problems with finite data the 
potential\index{potential} $\Phi$ has the form
given in \eqref{eq:potf2}
where $y\in\bbR^m$ is the data, $\cG:X \to \bbR^m$
is the forward model 
and $|\cdot|_{\Gamma}$ is
a covariance weighted norm on $\bbR^m$.
In this context the following corollary is useful.

\begin{corollary} \label{c:wp2}
Assume that $\Phi$ is given by \eqref{eq:potf2} and
that $\cG$ is approximated by a function
$\cG^N$ with the property that, for any $\epsilon>0$, 
there is $K'=K'(\epsilon)>0$ such that
\begin{equation}
\label{eq:forwardg}
|\G(u)-\G^N(u)|\le K'\exp\bigl(\epsilon\|u\nx^2\bigr)\psi(N)
\end{equation}
where $\psi(N) \to 0$ as $N \to \infty$. 
If $\cG$ and $\cG^N$ satisfy Assumptions \ref{ass:g1}
uniformly in $N$ then $\Phi$ and 
$\Phi^N:=\frac12|y-\cG^N(u)|_{\Gamma}^2$
satisfy the conditions necessary for application
of Theorem \ref{t:wp2} and all the conclusions of that
theorem apply.
\end{corollary}

\begin{proof} That (i), (ii) of Assumptions \ref{ass:1}
hold follows as in the proof of Lemma \ref{l:qwp}.
Also \eqref{eq:forward} holds since (for some $K(\cdot)$
defined in the course of the following chain of inequalities)
\begin{align*}
|\Phi(u)-\Phi^N(u)| & \le \frac12|2y-\cG(u)-\cG^N(u)|_{\Gamma}
|\cG(u)-\cG^N(u)|_{\Gamma}\\
& \le \Bigl(|y|+\exp\bigl(\epsilon\|u\nx^2+M\bigr)\Bigr) \times 
K'(\epsilon)\exp\bigl(\epsilon\|u\nx^2\bigr)\psi(N)\\
&\le K(2\epsilon)\exp(2\epsilon\|u\nx^2)\psi(N)
\end{align*}
as required.
\end{proof}

A notable fact concerning Theorem \ref{t:wp2} is that
the rate of convergence attained in the solution of the
forward problem, encapsulated in approximation of
the function $\Phi$ by $\Phi^N$,
is transferred into the rate of convergence
of the related inverse problem for measure $\mu$ given
by \eqref{eq:csagain2} and its approximation by $\mu^N.$
Key to achieving this transfer of rates of convergence
is the dependence of the constant in the forward error
bound \eqref{eq:forward} on $u$. In particular it is
 necessary that this constant is integrable by use of
the Fernique Theorem \ref{thm:fern}. In some
applications it is not possible to obtain such dependence.
Then convergence results can sometimes still be obtained, but
at weaker rates. We now describe a theory for this situation.

\begin{theorem}
\label{t:wp2f}
Assume that $\Phi$ and $\Phi^N$ satisfy 
Assumptions \ref{ass:1}(i),(ii) with constants uniform in $N$.
Assume also that, for any $R>0$ there is $K=K(R)>0$
such that, for all $u$ with $\|u\nx \le R$,
\begin{equation}
\label{eq:forwardf}
|\Phi(u)-\Phi^N(u)|\le K\psi(N)
\end{equation}
where $\psi(N) \to 0$ as $N \to \infty$. Then
the measures $\mu$ and $\mu^N$ are close 
with respect to the Hellinger distance: 
\begin{equation}
\label{eq:errorf}
\dhh(\mu,\mu^{N}) \to 0 
\end{equation}
as $N \to \infty.$
Consequently all moments of $\|u\nx$ under $\mu^N$ converge
to corresponding moments under $\mu$ as $N \to \infty.$ 
In particular the mean and, in the case $X$ is a Hilbert
space, the covariance operator, converge. 
\end{theorem}

\begin{proof}
Throughout the proof, all integrals are over $X$ unless
specified otherwise.
The constant $C$ changes from occurrence to occurrence.
The normalization constants $Z$ and $Z^N$ satisfy
lower bounds which are identical to that proved
for $Z$ in the course of establishing Theorem \ref{t:wp2}.

From Assumptions \ref{ass:1}(i) and \eqref{eq:forwardf}, 
\begin{align*}
|Z-Z^N| & \le \int_{X} |\exp\bigl(-\Phi(u)\bigr)-\exp\bigl(
-\Phi^N(u)\bigr)|d\mu_0\\
 &\le \int_{\{\|u\nx \le R\}} \exp\bigl(\epsilon\|u\nx^2-M\bigr)
|\Phi(u)-\Phi^N(u)|d\mu_0(u)\\
&\quad\quad
 +\int_{\{\|u\nx > R\}} 2\exp\bigl(\epsilon\|u\nx^2-M\bigr)d\mu_0(u)\\
& \le  \exp\bigl(\epsilon R^2-M\bigr) K(R)\psi(N)+J_{R}\\
& :=K_1(R)\psi(N)+J_{R}.
\end{align*}
Here
$$J_{R}=\int_{\{\|u\nx > R\}} 2\exp\bigl(\epsilon\|u\nx^2-M\bigr)d\mu_0(u).$$
Now, again by the Fernique Theorem \ref{thm:fern},
$J_{R} \to 0$ as $R \to \infty$ so, for any $\delta>0$,
we may choose $R>0$ such that $J_{R}<\delta.$ Now choose
$N>0$ so that $K_1(R)\psi(N)<\delta$ to deduce that
$|Z-Z^N|<2\delta.$ Since $\delta>0$ is arbitrary this proves
that $Z^N \to Z$ as $N \to \infty.$

From the definition of Hellinger distance we have
\begin{align*}
2\dhh(\mu,\mu^N)^2 & = \int\Bigl(
Z^{-\frac12}\exp\bigl(-\frac{1}{2}\Phi(u)\bigr)-(Z^N)^{-\frac12}
\exp\bigl(-\frac{1}{2} \Phi^N(u)\bigr)\Bigr)^2d\mu_0(u)\\
&\le I_1+I_2 
\end{align*}
where
\begin{align*}
I_1&=\frac{2}{Z}\int\Bigl(
\exp\bigl(-\frac{1}{2}\Phi(u)\bigr)-\exp(-\frac{1}{2}\Phi^N(u)\bigr)\Bigr)^2
d\mu_0(u),
\\
I_2&=2|Z^{-\frac12}-(Z^N)^{-\frac12}|^2\int \exp(-\Phi^N(u)\bigr)
d\mu_0(u).
\end{align*}

Now, again using Assumptions \ref{ass:1}(i) and 
equation \eqref{eq:forwardf}, 
\begin{align*}
I_1 & \le 
\frac{1}{2Z} \int_{\{\|u\nx \le R\}} K(R)^2\psi(N)^2\exp\bigl(\epsilon\|u\nx^2-M\bigr)d\mu_0(u)\\
&\quad\quad
+\frac{4}{Z} \int_{\{\|u\nx > R\}} 2\exp\bigl(\epsilon\|u\nx^2-M\bigr)d\mu_0(u)\\
&\le \frac{1}{2Z}K_2(R)\psi(N)^2+\frac{4}{Z}J_R, 
\end{align*}
for suitably chosen $K_2=K_2(R).$
An argument similar to the one above for $|Z-Z^N|$
shows that $I_1 \to 0$ as $N \to \infty.$

Note that the bounds on $Z,Z^N$ from below are
independent of $N$. Furthermore,
$$\int \exp\bigl(-\Phi^N(u)\bigr)d\mu_0(u) \le \int \exp\bigl(
\epsilon\|u\nx^2-M\bigr)d\mu_0(u)$$
with bound independent of $N$, by the Fernique Theorem \ref{thm:fern}.
Thus
$$|Z^{-\frac12}-(Z^N)^{-\frac12}|^2  \le C\bigl(Z^{-3} \vee (Z^N)^{-3}\bigr)|Z-Z^N|^2$$
and so $I_2 \to 0$ as $N \to \infty.$
Combining gives the desired continuity result in the
Hellinger metric.

The proof may be completed by the same arguments used
in Theorem \ref{t:wp2}. 
\end{proof}

\section{The Heat Equation}
\label{sec:heat}

Here we consider
a problem where the solution of the heat equation
is noisily observed at some fixed positive time $T>0.$
To be concrete we consider the heat equation
on a bounded open set $\Omega \subset
\bbR^d$, with Dirichlet boundary conditions, and 
written as an ODE in Hilbert space $\cH=L^2(\Omega)$:
\begin{equation}
\label{eq:heatII}
\frac{dv}{dt}+A v=0, \quad v(0)=u.
\end{equation}
Here 
$A=-\triangle$ with $D(A)=H^1_0(\Omega) \bigcap H^2(\Omega).$ 
We define the Sobolev spaces $\cHs$ as in \eqref{eq:Hs} with
$\cH=\cH^{0}=L^2(\Omega).$ 
We assume sufficient regularity conditions on $D$ and its 
boundary $\partial\Omega$ to ensure that the operator $A$
is the generator of an analytic semigroup
and we use \eqref{t:smooth} without comment
in what follows.

Assume that we observe the solution $v$ at time $T$, 
subject to error in the form of a Gaussian random field, 
and that we wish to recover the initial condition $u$. 
This problem is classically ill-posed, because the heat
equation is smoothing, and inversion of this
operator is not continuous on any Sobolev space $\cHs$. 
Nonetheless,
we will construct a well-defined Bayesian inverse problem.
We state a theorem showing
that the posterior measure is equivalent (in the
sense of measures) to the prior measure 
and then study the effect of approximation via a spectral method
in Theorem \ref{t:heaterr},
showing that the approximation error in the inverse problem
is exponentially small.

We place a prior measure on $u$ which is a Gaussian measure 
$\mu_0 \sim \cN(m_0,\cC_0)$ with $\cC_0=\beta A^{-\alpha},$ for
some $\beta>0,\alpha>\frac{d}{2}$. 
The lower bound on $\alpha$ ensures that samples
from the prior are continuous functions (Lemma \ref{lem:greg2}).

We assume that the observation is a function $y \in \cH$
and we model it as
\begin{equation}
\label{eq:hobsII}
y=e^{-A T}u +\eta
\end{equation}
where $\eta \sim \cN(0,\cC_1)$
and $\cC_1=\delta A^{-\gamma}$ for some $\delta>0$
and $\gamma>d/2$ so that $\eta$ is almost surely
continuous, by Lemma \ref{lem:greg2}. 
The forward model 
$\cG:\cH \to\cH$ is given by $\cG(u)= e^{-AT}u.$

By conditioning the Gaussian random variable $(u,y) \in
\cH \times \cH$ we find that the posterior measure for $u|y$ 
is also Gaussian $\cN(m,\cC)$ with mean
\begin{equation}
m=m_0+\frac{\beta}{\delta}e^{-AT}A^{\gamma-\alpha}
\Bigl(I+\frac{\beta}{\delta}e^{-2AT}A^{\gamma-\alpha}\Bigr)^{-1}
(y-e^{-AT}m_0)
\label{eq:fog}
\end{equation}
and covariance operator
\begin{equation}
\cC=\beta A^{-\alpha}
\Bigl(I+\frac{\beta}{\delta}e^{-2AT}A^{\gamma-\alpha}\Bigr)^{-1}.
\label{eq:hevII}
\end{equation}

We can also derive a formula for
the Radon-Nikodym derivative between
$\muy(du)=\bbP(du|y)$ and the prior $\mu_0(du).$
We define $\Phi:\cH \times\cH \to \bbR$ by 
\begin{equation}
\Phi(u;y)=
\frac{1}{2}\|\cC_1^{-\frac12}e^{-AT}u\|^2-\langle \cC^{-\frac12}_1 e^{-AT}u,\cC^{-\frac12}_1 y\rangle.
\label{eq:bheatII}
\end{equation}
It is a straightforward application
of the theory of Gaussian measures \cite{bog98,DapZab92}, 
using the continuity properties 
of $\Phi$ established below, to prove the following:

\begin{theorem} \cite{Stuart10}
Consider the inverse problem for the initial
condition $u$ in \eqref{eq:heatII},
subject to observation in the form \eqref{eq:hobsII} with
observational noise $\eta \sim \cN(0,\delta A^{-\gamma})$,
$\delta>0$ and $\gamma>\frac{d}{2}.$ Assume that the 
prior measure is a Gaussian $\mu_0=\cN(m_0,\beta A^{-\alpha})$ 
with $m_0 \in \cH^{\alpha},\beta>0$ and $\alpha>\frac{d}{2}$.
Then the posterior measure $\muy$ is Gaussian
with mean and variance determined by  
\eqref{eq:fog} and \eqref{eq:hevII}. Furthermore, $\muy$
and the prior measure $\mu_0$ are equivalent Gaussian measures
with Radon-Nikodym derivative \eqref{eq:cs} given by \eqref{eq:bheatII}.
\end{theorem}

Now we study the properties of $\Phi.$
To this end it is helpful to define, for any
$\llambda>0$, the compact
operator $K_{\llambda}:\cH \to \cH$ given by 
$$K_{\llambda}:=\cC_1^{-\frac12}e^{-\llambda AT}.$$ 
Note that, for any $0<\theta_1<\theta_2<\infty$ there is $C>0$
such that, for all $u \in \cH$,
$$\|K_{\theta_2}u\| \le C \|K_{\theta_1}u\|.$$ 

\begin{lemma} \label{l:hIg} The function $\Phi:\cH\times\cH \to \bbR$
satisfies Assumptions \ref{ass:1} with $X=Y=\cH$ and, furthermore, 
for any $\epsilon \in (0,1)$, there is $C=C(\epsilon)$ such that
$$|\Phi(u;y)-\Phi(v;y)| \le C\Bigl(\|K_1 u\|+\|K_1 v\|+
\|K_{\epsilon} y\|
\Bigr)\|K_{1-\epsilon}u-K_{1-\epsilon}v\|.$$
In particular,
$\Phi(\cdot;y):\cH \to \bbR$
is continuous. 
\end{lemma}

\begin{proof} We may write
$$\Phi(u;y)=
\frac{1}{2}\|\cC_1^{-\frac12}e^{-AT}u\|^2-
\langle \cC_1^{-\frac12}e^{-\frac12 AT}u,
\cC_1^{-\frac12}e^{-\frac12 AT}y\rangle.$$
By the Cauchy-Schwarz inequality we have, for any $\delta>0$,
$$\Phi(u;y) \ge -\frac{\delta^2}{2}
\|\cC_1^{-\frac12}e^{-\frac12 AT}u\|^2-\frac{1}{2\delta^2}
\|\cC_1^{-\frac12}e^{-\frac12 AT}y\|^2
$$
so that, by the compactness of $K_{\frac12}$, 
Assumption \ref{ass:1}(i) holds.
Assumption \ref{ass:1}(ii) holds, by a similar
Cauchy-Schwarz argument, with
$$\Phi(u;y) \le \frac{1}{2}
\|\cC_1^{-\frac12}e^{-AT}u\|^2+
\frac{1}{2}\|\cC_1^{-\frac12}e^{-\frac12 AT}y\|^2
+
\frac{1}{2}\|\cC_1^{-\frac12}e^{-\frac12 AT}u\|^2$$
so that, by the compactness of $K_{\llambda}$, 
\begin{equation}
\Phi(u;y) \le C\Bigl(1+\|u\|^2\Bigr).
\label{eq:see}
\end{equation}

Note that
$$\langle \cC_1^{-\frac12}e^{-\frac12 AT}u,
\cC_1^{-\frac12}e^{-\frac12 AT}y\rangle
=\langle \cC_1^{-\frac12}e^{-(1-\epsilon) AT}u,
\cC_1^{-\frac12}e^{-\epsilon AT}y\rangle.$$
Since $\Phi$ is quadratic in $u$
the desired Lipschitz property holds.
\end{proof}

Now we consider approximation of the posterior measure $\muy$
given by \eqref{eq:bheatII}. Specifically we define
$P^N$ to be orthogonal projection in $\cH$
into the subspace $\{\phi_k\}_{|k|\le N}$ (a subset of
the eigenfunctions of $A$ as defined just before \eqref{eq:carf2}) 
and define the measure $\mu^N$ given by 
\begin{equation}
\frac{d\mu^N}{d\mu_0}(u) \propto \exp \Bigl( -\Phi(P^N u;y)\Bigr).
\label{eq:bheatIII}
\end{equation}
The measure $\mu^N$ is identical to $\mu_0$ on the orthogonal
complement of $P^N\cH$. 
We now use the theory from the preceding section to
estimate the distance between $\muy$ and $\mu^N$.

\begin{theorem} There are constants $c_1>0, c_2>0$,
independent of $N$, such that $\dhh(\muy,\mu^N)
\le c_1 \exp(-c_2 N^2).$ Consequently the mean and
covariance operator of $\muy$ and $\mu^N$ are
${\cal O}\bigl(\exp(-c_2 N^2)\bigr)$ close in
the $\cH$ and $\cH-$operator norms respectively.
\label{t:heaterr}
\end{theorem}

\begin{proof} We apply Theorem \ref{t:wp2} with $X=\cH$.
By Lemma \ref{l:hIg}, 
together with the fact that $\|P^N u\| \le \|u\|$, we
deduce that Assumptions \ref{ass:1} hold for $\Phi$ and
$\Phi^N$, with constants independent of $N$. Furthermore,
from the Lipschitz bound in Lemma \ref{l:hIg}, we have
$$|\Phi(u;y)-\Phi^N(u;y)| \le C\Bigl(\|u\|+\|y\|\Bigr)
\|K_{\frac12}(u-P^N u)\|.$$
But
$$\|K_{\frac12}(u-P^N u)\|^2=\frac{1}{\delta}\sum_{|k|>N}  
\lambda_{k}^{\gamma} \exp(-\lambda_k T)|u_k|^2.$$
Since the eigenvalues $\lambda_k$ grow like $|k|^2,$
and since $x^{\gamma}\exp(-x T)$ is monotonic decreasing
for $x$ sufficiently large, we deduce that 
$$\|K_{\frac12}(u-P^N u)\|^2 \le 
c_1 \exp(-c_2 N^2)
\sum_{|k|>N}|u_k|^2 \le
c_1 \exp(-c_2 N^2)\|u\|^2.$$ 
The result follows (possibly by redefinition of $c_1, c_2$).
\end{proof}

\section{Lagrangian Data Assimilation}
\label{sec:lag}

In this section we turn to a non-Gaussian
nonlinear example where the full power of the abstract theory
is required.
In oceanography a commonly used method of gathering data about
ocean currents, temperature, salinity and so forth
is through the use of Lagrangian instruments:
objects transported by the fluid velocity field, 
which transmit positional information using GPS. 
The inverse problem termed {\em Lagrangian data
assimilation} is to determine
the velocity field in the ocean from the Lagrangian 
data \cite{IKJ02,KIJ03}.

In this section we study an idealized model which captures
the essence of Lagrangian data assimilation
as practised in oceanography. For the fluid
flow model we use the {\em Stokes equations}\index{Stokes
equations}, describing incompressible Newtonian fluids
at moderate Reynolds number. The real equations of oceanography
are, of course, far more complex, requiring evolution of coupled fields
for velocity, temperature and salinity. However the dissipative
and incompressible nature of the flow field for the Stokes
equations captures the key mathematical properties of ocean flows,
and hence provides a useful simplified model.

We consider the incompressible
Stokes equations written in the form:
\begin{subequations}
\label{eq:stokes11}
\begin{equation}
\label{eq:stokesa}
\frac{\partial v}{\partial t}=\nu \Delta v-\nabla p+f,
\quad(x,t) \in \Omega\times [0,\infty),
\end{equation}
\begin{equation}
\nabla \cdot v=0, \quad (x,t) \in \Omega \times [0,\infty),
\label{eq:stokesb}
\end{equation}
\begin{equation}
v=u, \quad (x,t) \in {\overline \Omega} \times \{0\}.
\label{eq:stokesc}
\end{equation}
\end{subequations}
Here $\Omega$ is the unit square.
We impose periodic boundary conditions on the velocity field $v$ and
the pressure $p$.
We assume that $f$ has zero average over $\Omega$, noting that
this implies the same for $v(x,t)$, provided 
that $u(x)=v(x,0)$ has zero initial average.
See \cite{temam,temam_ns} for
definitions of
the Leray projector $P:L^2_{per}\to \cH$
and Stokes operator $A$. 
We employ the Hilbert spaces $\{\cH^s,\|\cdot\|_s\}$
defined by \eqref{eq:Hs} 
and note that $\cH^s=D(A^{s/2})$ 
for any $s>0.$

The PDE can be formulated as a linear dynamical system
on the Hilbert space 
\begin{equation}
\label{eq:spaceh}
\cH = \Bigl\{ u \in
L^2_{\rm{per}}(\Omega)\Big|\int_{\Omega}u dx = 0, \, \nabla \cdot u
= 0\Bigr\},
\end{equation}
with the usual $L^2(\Omega)$
norm and inner-product on this subspace of $L^2_{\rm per}(\Omega).$
If we let $\psi=Pf$ then we may write the 
equation \eqref{eq:stokes11} 
as an ODE in Hilbert space $\cH:$
\begin{equation}
\frac{dv}{dt}+\nu Av=\psi, \quad v(0)=u.
\label{eq:stokes1}
\end{equation}
We assume that we are
given noisy observations of $J$ Lagrangian 
tracers with positions $z_j$ solving
the integral equations
\begin{equation}
z_j(t)=z_{j,0}+\int_0^t v(z_j(s),s)ds.
\label{eq:tracers2}
\end{equation}

For simplicity assume that we observe all the tracers $z$ at the
same  set of positive times $\{t_k\}_{k=1}^K$ and that the initial
particle tracer positions $z_{j,0}$ are known to us:
\begin{equation}
\label{eq:ldata}
y_{j,k}=z_j(t_k)+\eta_{j,k}, \quad j=1,\dots, \jb\,\,\,k=1, \dots, \kb,
\end{equation}
where the $\eta_{j,k}$'s are zero mean Gaussian random variables.
Concatenating data we may write
\begin{equation}
\label{eq:ldata2}
y={\cal G}(u)+\eta
\end{equation}
with $y=(y_{1,1}, \dots, y_{\jj,\kk})^*$ and 
$\eta \sim {\cal N}(0,\Gamma)$ for some covariance matrix $\Gamma$
capturing the correlations present in the noise.
Note that ${\cal G}$ is a complicated function
of the initial condition for the Stokes equations,
describing the mapping from this initial condition
into the positions of Lagrangian trajectories at
positive times.
We will show that the function ${\cal G}$ maps of
$\cH$ into $\bbR^{2\jj\kk}$, and is continuous on
a dense subspace of $\cH$.

The objective of the inverse problem 
is to find the initial velocity field $u$, given $y$.
We adopt a Bayesian approach and identify
$\mu(du)=\bbP(u|y)du.$
We now spend some time developing the Bayesian framework,
culminating in Theorem \ref{thm:ysae} which shows that $\mu$
is well-defined. The reader interested purely in approximation
of $\mu$ can skip straight to Theorem \ref{t:stokerr}.

The following result shows that the tracer
equations \eqref{eq:tracers2} have a solution,
under mild regularity assumptions on the
initial data.
An analogous result is proved in \cite{DR07} for the
case where the velocity field is governed by
the Navier-Stokes equation and the proof may be easily
extended to the case of the Stokes equations.

\begin{theorem} \label{t:trace}
Let $\psi\in L^2(0,T;\cH)$
and let $v\in C([0,T];\cH)$ denote the solution of
\eqref{eq:stokes1} with initial data $u \in \cH$.
Then the integral equation
\eqref{eq:tracers2}
has a unique solution $z \in C([0,T],\bbR^2).$
\end{theorem} 

We assume throughout that $\psi$
is sufficiently regular that this
theorem applies.
To determine a formula for the 
probability of $u$ given $y$, we apply the Bayesian
approach described in \cite{CDRS08} for the Navier-Stokes
equations, and easily generalized to the Stokes equations.
For the prior measure we take 
$\mu_0=\cN(0,\beta A^{-\alpha})$ for some $\beta>0,\alpha>1$,
with the condition on $\alpha$ chosen to ensure
that draws from the prior are in $\cH$, by Lemma \ref{lem:greg2}.
We condition the prior on the
observations, to find the {\em posterior} measure on $u$.
The likelihood of $y$ given $u$ is
$$\bbP\,(y\mid u) \propto \exp \Bigl(-\frac12|y-{\cal G}(u)|_{\Gamma}^2\Bigr).$$
This suggests the formula
\begin{equation}
\label{eq:RNlda}
\frac{d\muy}{d\mu_0}(u) \propto \exp\Bigl(-\Phi(u;y)\Bigr)
\end{equation}
where
\begin{equation}
\label{eq:lphi}
\Phi(u;y):=\frac12|y-{\cal G}(u)|_{\Gamma}^2
\end{equation}
and $\mu_0$ is the prior Gaussian measure. 
We now make this assertion rigorous. The first step is to
study the properties of the forward model $\cG.$
Proof of the following lemma is given after statement
and proof of the main approximation result, Theorem \ref{t:stokerr}.

\begin{lemma} \label{l:5a} Assume that
$\psi \in C([0,T];\cH^{\gamma})$ for some $\gamma \ge 0.$
Consider the forward model 
$\cG:\cH \to \bbR^{2\jj\kk}$ defined 
by \eqref{eq:ldata},\eqref{eq:ldata2}. 
\begin{itemize}
\item If $\gamma \ge 0$ then, for any $\ell \ge 0$ 
there is $C>0$ such that,
for all $u \in \cH^{\ell},$
$$|\cG(u)| \le C\bigl(1+\|u\|_{\ell}\bigr).$$
\item If $\gamma>0$ then, for any $\ell>0$ and $R>0$ 
and for all $u_1,u_2$ with
$\|u_1\|_{\ell} \vee \|u_2\|_{\ell}<R$,
there is $L=L(R)>0$ such that
$$|\cG(u_1)-\cG(u_2)| \le L\|u_1-u_2\|_{\ell}.$$
Furthermore, for any $\epsilon>0$, there is $M>0$ such that
$L(R) \le M\exp(\epsilon R^2).$
\end{itemize}
Thus $\cG$ satisfies Assumptions \ref{ass:g1} with
$X=\cH^{s}$ and any $s\ge 0.$
\end{lemma}

Since $\cG$ is continuous on $\cH^{\ell}$
for $\ell>0$ and since, by Lemma \ref{lem:greg2},
draws from $\mu_0$ are almost surely in $\cH^s$ for
any $s<\alpha-1,$ use of the
techniques in \cite{CDRS08}, employing the Stokes
equation in place of the Navier-Stokes equation,
shows the following:

\begin{theorem} \label{thm:ysae} Assume that
$\psi \in C([0,T];\cH^{\gamma})$, for some $\gamma>0$,
and that the
prior measure $\mu_0=\cN(0,\beta \cA^{-\alpha})$ is
chosen with  $\beta>0$ and $\alpha>1.$
Then the measure $\muy(du)=\bbP(du|y)$
is absolutely continuous with respect to the prior $\mu_0(du)$,
with Radon-Nikodym derivative given by \eqref{eq:RNlda}.
\end{theorem}

In fact the theory in \cite{CDRS08} may be used to
show that the measure $\mu$ is Lipschitz in the data
$y$, in the Hellinger metric. This well-posedness
underlies the following study of the 
approximation of $\muy$ in a finite dimensional space. 
We define $P^N$ to be orthogonal projection in $\cH$
into the subspace $\{\phi_k\}_{|k|\le N}$; recall
that $k \in \bbK:=\bbZ^2\backslash\{0\}.$
Since $P^N$ is an orthogonal projection
in any $\cH^a$ we have $\|P^N u\nx \le \|u\nx.$  
Define
$$\cG^N(u):=\cG(P^N u).$$
The approximate posterior measure $\mu^N$ is given
by \eqref{eq:RNlda} with $\cG$ replaced by $\cG^N$.
As in the last section it is identical to the prior
on the orthogonal complement of $P^N \cH$. On $P^N\cH$
itself the measure is finite dimensional and amenable to
sampling techniques as demonstrated in \cite{CDRS09}.
We now quantify the error arising from approximation
of $\cG$ in the finite dimensional subspace $P^N X.$

\begin{theorem} Let the assumptions of Theorem \ref{thm:ysae}
hold. Then, for any $q<\alpha-1$, there
is a constant $c>0$, independent of $N$, such that 
$\dhh(\muy,\mu^N) \le c N^{-q}.$ 
Consequently the mean and
covariance operator of $\muy$ and $\mu^N$ are
${\cal O}\bigl(N^{-q}\bigr)$ 
close in the $\cH$ and $\cH-$operator norms respectively.
\label{t:stokerr}
\end{theorem}

\begin{proof} 
We set $X=\cH^s$ for any $s \in (0,\alpha-1).$
We employ Corollary \ref{c:wp2}.
Clearly, since $\cG$ satisfies Assumptions \ref{ass:g1}
by Lemma \ref{l:5a}, so too does $\cG^N$, with constants
uniform in $N.$
%We set $X'=\cH^{\ell}$ for $\ell \in (0,s).$
It remains to establish \eqref{eq:forwardg}.
Write $u \in \cH^s$ as
$$u=\sum_{k \in \bbK} u_k \phi_k$$
and note that
$$\sum_{k \in \bbK} |k|^{2s}|u_k|^2<\infty.$$ 
We have, for any $\ell \in (0,s)$,
\begin{align*}
\|u-\pn u\|_{\ell}^2 &=\sum_{|k|>N}|k|^{2\ell}|u_k|^2\\
&=\sum_{|k|>N}|k|^{2(\ell-s)}|k|^{2s}|u_k|^2\\
&\le N^{-2(s-\ell)}\sum_{|k|>N}|k|^{2s}|u_k|^2\\
&\le C\|u\|_{s}^2 N^{-2(s-\ell)}.
\end{align*}
By the Lipschitz properties of $\cG$ from Lemma \ref{l:5a}
we deduce that, for any $\ell \in (0,s)$, 
\begin{align*}
|\cG(u)-\cG(P^N u)| & \le M\exp\bigl(\epsilon\|u\|_{\ell}^2\bigr)
\|u-P^Nu\|_{\ell}\\
 & \le C^{\frac12}M\exp\bigl(\epsilon\|u\|_{s}^2\bigr)\|u\|_{s}
N^{-(s-\ell)}.
\end{align*}
This establishes the desired error bound \eqref{eq:forwardg}.
It follows from Corollary \ref{c:wp2}
that $\mu^N$ is ${\cal O}\bigl(N^{-(s-\ell)}\bigr)$ 
close to $\mu$ in the Hellinger distance. 
Choosing $s$ arbitrarily close to its upper bound, 
and $\ell$ arbitrarily close to zero, yields the
optimal exponent $q$ as appears in the theorem statement.
\end{proof}

\begin{proof} {\em of Lemma \ref{l:5a}}
Throughout the proof, the constant $C$ may change
from instance to instance, but is always independent of the $u_i.$
It suffices to consider a single observation
so that $J=K=1.$ Let $\zi(t)$ solve
$$\zi(t)=\zi_0+\int_0^t \vi(\zi(\tau),\tau)d\tau$$
where $\vi(x,t)$ solves \eqref{eq:stokes11} with $u=u_i.$ 

Let $\ell \in [0,2+\gamma).$
Recall that, by \eqref{t:thirteen},
\begin{equation}
\label{eq:buv}
\|\vi(t)\|_s \le C\Bigl(\frac{1}{t^{(s-\ell)/2}}\|u_i\|_{\ell}
+\|\psi\|_{C([0,T];\cH^{\gamma})}\Bigr),
\end{equation}
for $s \in [\ell,2+\gamma).$
Also, by linearity and \eqref{t:smooth},
\begin{equation}
\label{eq:buv2}
\|v^{(1)}(t)-v^{(2)}(t)\|_s \le \frac{C}{t^{(s-\ell)/2}}\|u_1-u_2\|_{\ell}.
\end{equation}

To prove the first part of the lemma note that, by the
Sobolev embedding Theorem, for any $s>1$,
\begin{align*}
|\zi(t)|&
\le|\zi_0|+\int_0^t\|\vi(\cdot,\tau)\|_{L^{\infty}}d\tau\\
&\le C\Bigl(1+\int_0^t\|\vi(\cdot,\tau)\|_{s}d\tau\Bigr)\\
&\le C\Bigl(1+\int_0^t\frac{1}{\tau^{(s-\ell)/2}}\|u_i\|_{\ell}d\tau\Bigr).
\end{align*}
For any $\gamma \ge 0$ and 
$\ell \in [0,2+\gamma)$ we may choose $s$ such that $s\in 
[\ell,2+\gamma)\bigcap (1,\ell+2).$
Thus the singularity is integrable and we have, for any $t\ge 0$,
$$|\zi(t)| \le C\bigl(1+\|u_i\|_{\ell}\bigr)$$
as required.

To prove the second part of the lemma choose $\ell \in (0,2+\gamma)$
and then choose $s\in[\ell-1,1+\gamma)\cap(1,\ell+1)$;
this requires $\gamma>0$ to ensure a nonempty 
intersection. Then
\begin{equation}
\label{eq:buv11}
\|\vi(t)\|_{1+s} \le C\Bigl(\frac{1}{t^{(1+s-\ell)/2}}\|u_i\|_{\ell}
+\|\psi\|_{C([0,T];\cH^{\gamma})}\Bigr).
\end{equation}
Now we have 
\begin{align*}
|\zo(t)-\zt(t)|&\le|\zo(0)-\zt(0)|+\int_0^t|\vo(\zo(\tau),\tau)-\vt(\zt(\tau),\tau)|d\tau\\
& \le\int_0^t\|D\vo(\cdot,\tau)\|_{L^\infty}|\zo(\tau)-
\zt(\tau)|d\tau\\
&\quad+\int_0^t\|\vo(\cdot,\tau)-\vt(\cdot,\tau)\|_{L^\infty}d\tau\\
& \le\int_0^t\|\vo(\cdot,\tau)\|_{1+s}|\zo(\tau)-
\zt(\tau)|d\tau\\
&\quad+\int_0^t\|\vo(\cdot,\tau)-\vt(\cdot,\tau)\|_{s}d\tau\\
& \le\int_0^t C\Bigl(\frac{1}{\tau^{(1+s-\ell)/2}}\|\uo\|_{\ell}+
\|\psi\|_{C([0,T];\cH^{\gamma})}\Bigr)|\zo(\tau)-\zt(\tau)|d\tau\\
&\quad+\int_0^t\frac{C}{\tau^{(s-\ell)/2}}\|\uo-\ut\|_{\ell}d\tau.
\end{align*}
Both time 
singularities are integrable and application of the
Gronwall inequality\index{Gronwall inequality} from
Lemma \ref{lem:gron} gives, for some $C$ depending on
$\|\uo\|_{\ell}$ and $\|\psi\|_{C([0,T];\cH^{\gamma})}$, 
\begin{equation*}
\|\zo-\zt\|_{L^{\infty}((0,T);\bbR^2)} \le C \|\uo-\ut\|_{\ell}.
\end{equation*}
The desired Lipschitz bound on $\cG$ follows.
In particular, the desired dependence of the Lipschitz
constant follows from the fact that, for any $\epsilon>0$
there is $M>0$ with the property that, for all $\llambda \ge 0$,
$$1+\llambda\exp(\llambda) \le M\exp(\eps\llambda^2).$$
\end{proof}

We conclude this section with the results of numerical experiments
illustrating the theory. We compute the posterior distribution on
the initial condition for Stokes equations from observation of
$J$ Lagrangian trajectories at one time $t=0.1$. 
The prior measure is taken to be $\cN(0,400\times A^{-2})$.
The initial condition used to generate the data is found
by making a single draw from the prior measure and
the observational noise on the Lagrangian data is i.i.d
$\cN(0,\gamma^2)$ with $\gamma=0.01$.

Note that, in the periodic geometry
assumed here, the Stokes equations can be solved exactly by
Fourier analysis \cite{temam_ns}. Thus there
are four sources of approximation when attempting to sample
the posterior measure on $u.$ These are
\begin{itemize}
\item (i) the effect of generating approximate samples from
the posterior measure by use of MCMC methods; 
\item (ii) the effect of approximating $u$ in a finite space
found by orthogonal projection on the eigenbasis of the
Stokes operator;
\item (iii) the effect of interpolating a velocity field
on a grid, found from use of the FFT,
into values at the arbitrary locations of Lagrangian tracers;
\item (iv) the effect of time-step in an Euler integration
of the Lagrangian trajectory equations.
\end{itemize}

The MCMC method that we use is a generalization of the random
walk Metropolis method and is detailed in \cite{CDRS09}.
The method is appropriate for sampling measures
absolutely continuous with respect to a Gaussian in
the situation where it is straightforward to sample
directly from the Gaussian itself. We control the
error (i) simply by running the MCMC method until time
averages of various test statistics have converged; the
reader interested in the effect of this Monte Carlo
error should consult \cite{CDRS09}. The
error in (ii) is precisely the error which we
quantify in Theorem \ref{t:stokerr}; for the particular
experiments used here we predict an error of order
$N^{-q}$ for any $q \in (0,1).$ In this paper we have not
analyzed the errors resulting from (iii) and (iv): 
these approximations are not included in the analysis
leading to Theorem \ref{t:stokerr}. However we anticipate
that Theorem \ref{t:wp2} or Theorem \ref{t:wp2f}
could be used to study such approximations and the
numerical evidence which follows below is consistent with
this conjecture.

\begin{figure}[htp]
\begin{center}
\scalebox{0.6}{\includegraphics{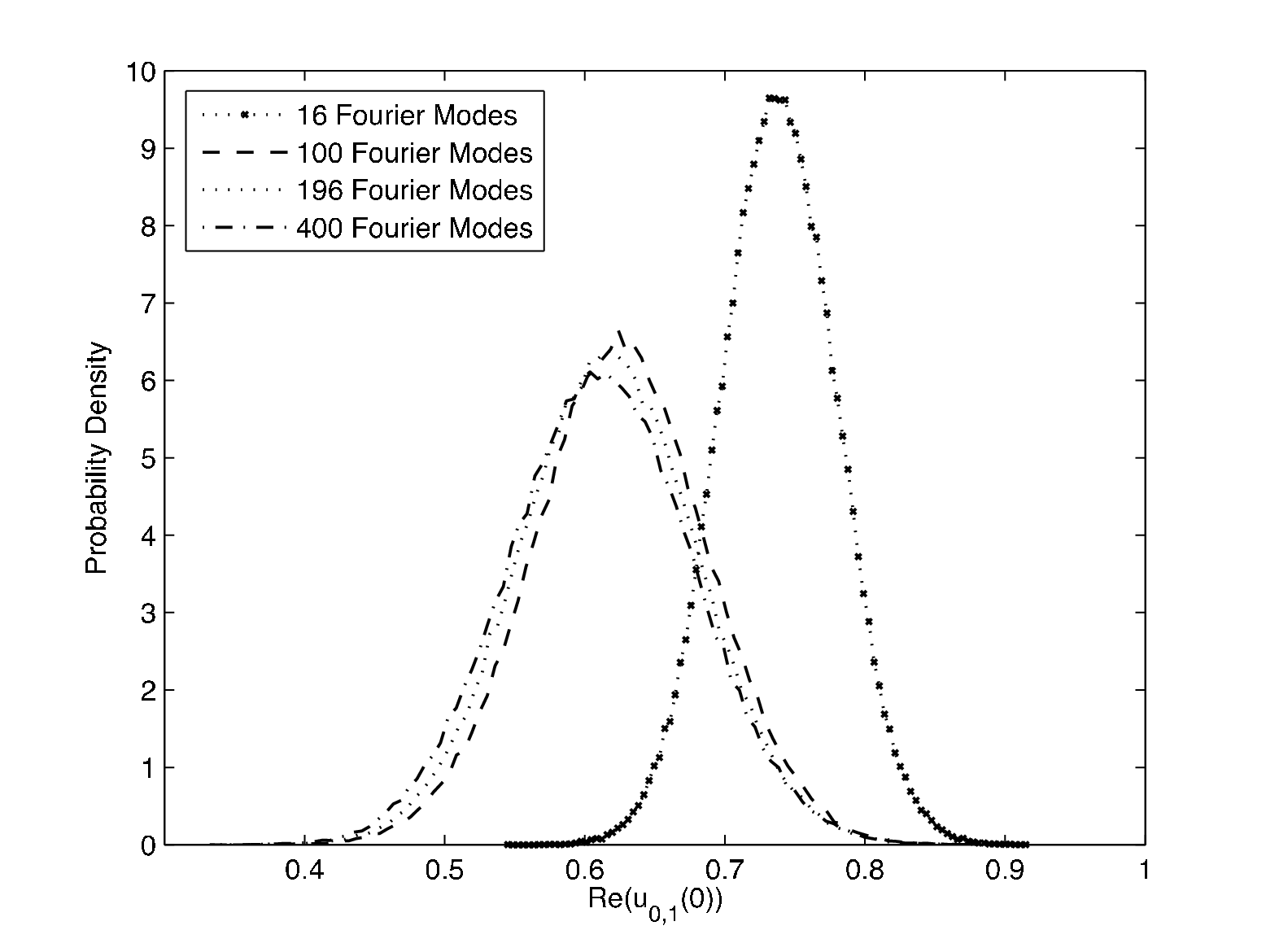}}
\end{center}
\caption{Marginal distributions on ${\rm{Re}}(u_{0,1}(0))$ with
  differing numbers of Fourier modes. \label{JLFR01}}
\end{figure}

\begin{figure}[htp]
\begin{center}
\scalebox{0.6}{\includegraphics{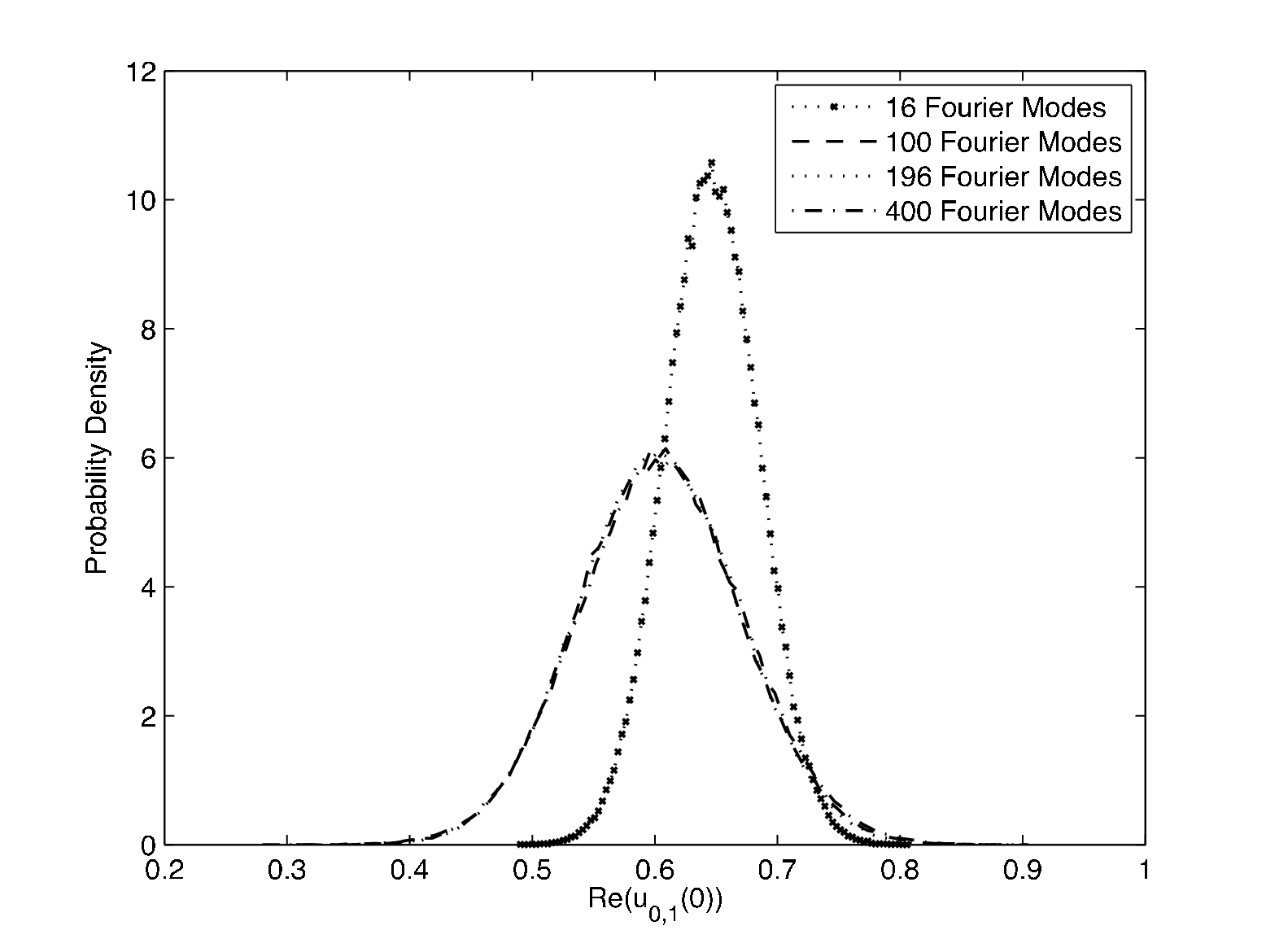}}
\end{center}
\caption{Marginal distributions on ${\rm{Re}}(u_{0,1}(0))$ with
  differing numbers of Fourier modes, bicubic
  interpolation used. \label{JLBCR01}}
\end{figure}

\begin{figure}[htp]
\begin{center}
\scalebox{0.6}{\includegraphics{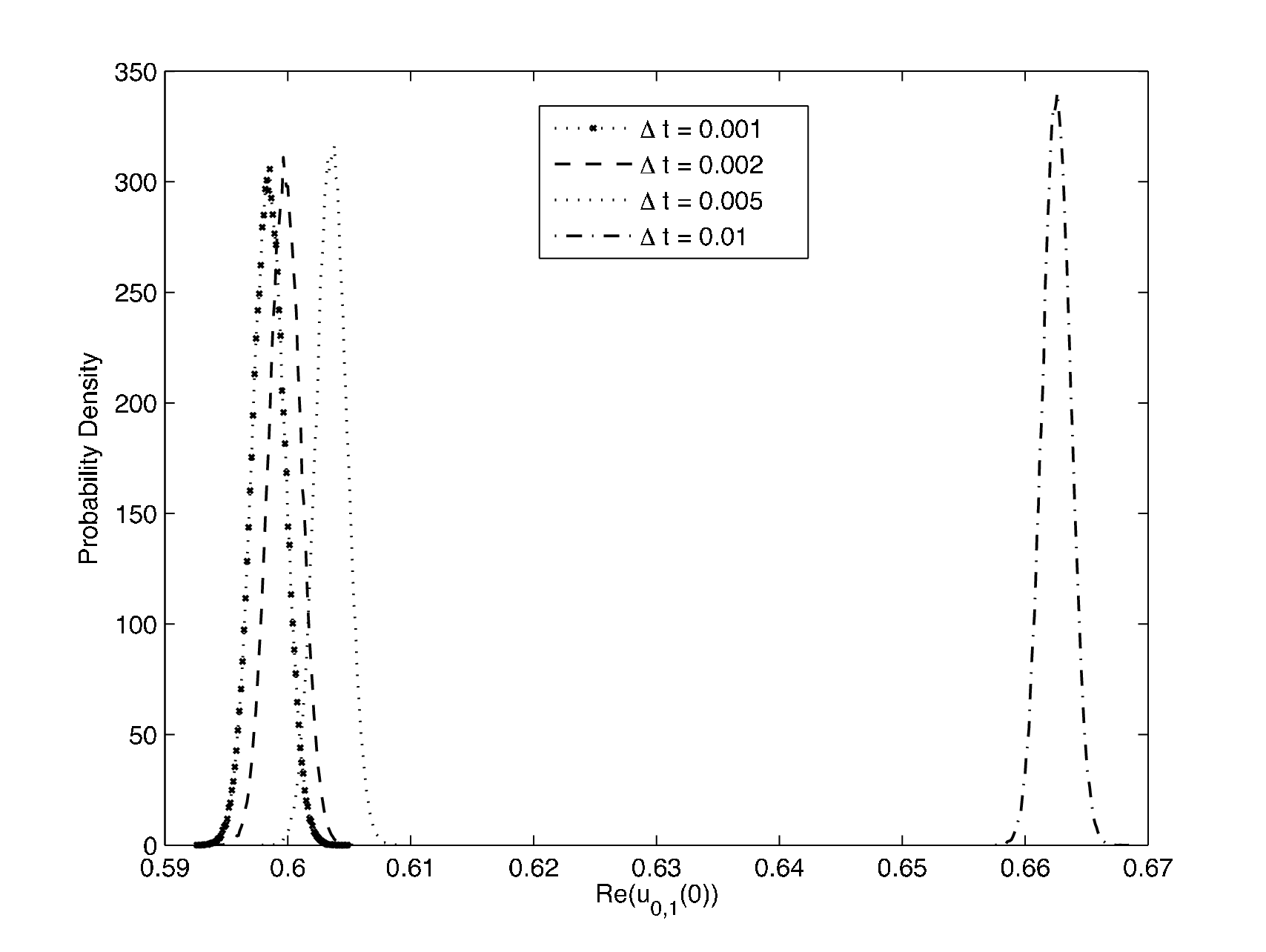}}
\end{center}
\caption{Marginal distributions on ${\rm{Re}}(u_{0,1}(0))$ with
  differing timestep, Lagrangian data \label{RefineDtLag}}
\end{figure}

In the following three numerical experiments (each
illustrated by a figure) we study the
effect of one or more of the approximations (ii), (iii) and
(iv) on the empirical distribution (`histogram') found
from marginalizing data from the MCMC method onto the
real part of the Fourier mode with wavevector $k=(0,1).$
Similar results are found for other Fourier modes although
it is important to note that
at high values of $|k|$ the data is uninformative and the
posterior is very close to the prior (see \cite{CDRS09}
for details).
The first two figures use $J=9$ Lagrangian trajectories,
whilst the third uses $J=400.$
Figure \ref{JLFR01} shows the effect of increasing the number of
Fourier modes\footnote{Here by number of Fourier modes, we mean the dimension of the Fourier space approximation, ie then number of grid points} used from $16$, through $100$ and $196$, to
a total of $400$ modes and illustrates Theorem \ref{t:stokerr}
in that convergence to a limit is observed as the number
of Fourier modes increases.
However this experiment is conducted by using bilinear
interpolation of the velocity field on the grid, in order
to obtain off-grid velocities required for
particle trajectories. At the
cost of quadrupling the number of FFTs it is possible
to implement bicubic interpolation \footnote{Bicubic interpolation with no added FFTs is also possible by using finite difference methods to find the partial derivatives, but at a lower order of accuracy}. Conducting the
same refinement of the number of Fourier modes then yields
Figure \ref{JLBCR01}. 
Comparison of Figures \ref{JLFR01} and \ref{JLBCR01} shows
that the approximation (iii) by increased order
of interpolation leads to improved approximation of 
the posterior distribution, and Figure \ref{JLBCR01} alone 
again illustrates Theorem \ref{t:stokerr}. 
Figure \ref{RefineDtLag} shows the effect (iv) of reducing the time-step
used in the integration of the Lagrangian trajectories.
Note that many more ($400$) particles were used to
generate the observations leading to this 
figure than were used in the preceding two figures. 
This explains
the quantitatively different
posterior distribution; in particular the variance 
in the posterior distribution is considerably smaller.
The result shows clearly that reducing the time-step
leads to convergence in the posterior distribution.

\section{Eulerian Data Assimilation}
\label{sec:eul}

In this section we consider a data assimilation problem that is 
related to weather forecasting applications.
In this problem, direct observations are made of the velocity field
of an incompressible viscous flow at some fixed points in 
space-time, 
the mathematical model is the two-dimensional Navier-Stokes
equations on a torus, 
and the objective is to obtain an estimate of the 
initial velocity field.
The spaces $\cH$ and $\cH^s$ are as defined in Section \ref{sec:lag},
with $\|\cdot\|_s$ the norm in $\cH^s$ and
$\|\cdot\|=\|\cdot\|_0$.
The definitions of $A$, the Stokes operator, and $P$,
the Leray projector, are also as in the 
previous section \cite{temam, temam_ns}. 

We consider the incompressible two-dimensional
Navier-Stokes equations
\begin{align*}
\frac{\partial v}{\partial t}=\nu \lap v-(v\cdot \grad)v-\grad p
      &+f,\quad (x,t)\in D\times[0,\infty), \\
\cdiv v&=0, \quad (x,t)\in D\times[0,\infty), \\
v&=u, \quad (x,t)\in \overline D\times \{0\},
\end{align*}
where $D$ is a unit square as before and the boundary conditions are periodic.
We apply the Leray Projector $P:L^2_{per}(\Omega)\to\cH$ and write
the Navier-Stokes equations as an ordinary differential equation in $\cH$
\begin{align}\label{eq:nsef}
\frac{\ud v}{\ud t}+\nu Av+B(v,v)=\psi,\quad v(0)=u
\end{align}
with $A$ the Stokes operator, $B(v,v)=P((v\cdot \grad) v)$ and 
$\psi=Pf$. 

For simplicity we assume that we make noisy observations of 
the velocity field $v$ at time $t>0$ and at points
$x_1,\dots,x_K\in\Omega$:
$$
y_k=v(x_k,t)+\eta_k,\quad k=1,\dots,K.
$$ 
We assume that the noise is Gaussian and the $\eta_k$ form
an i.i.d sequence with $\eta_1\sim \cN(0,\gamma^2)$. 
It is known (see Chapter 3 of \cite{temam}, for example) that for $u\in \cH$ and $f\in L^2(0,T;\cH^s)$ with $s>0$ a unique solution 
to (\ref{eq:nsef}) exists which satisfies $u\in L^\infty(0,T;\cH^{1+s})\subset L^\infty(0,T;L^\infty(D))$. 
Therefore for such initial condition and forcing function the value of $v$ at any $x\in\Omega$ 
can be written as a function of $u$.
Hence, we can write
$$
y=\cG(u)+\eta
$$
where  $y=(y_1,\cdots,y_{\kk})^T\in\bbR^{\kk}$ and $\eta=(\eta_1,\dots, \eta_k)^T\in\bbR^{\kk}$
is distributed as $\cN(0,\gamma^2 I)$ and
\begin{equation}
\label{eq:bvpa}
\cG(u)=(v(x_1,t),\cdots,v(x_{\kk},t))^T.
\end{equation}
Now consider a Gaussian prior measure $\mu_0\sim \cN(u_b,\beta A^{-\alpha})$
with $\beta>0$ and $\alpha>1$; recall that the second condition
ensures that functions drawn from the prior are in $\cH$, by
Lemma \ref{lem:greg2}.
In Theorem 3.4 of \cite{CDRS08} 
it is shown that with such prior measure, the 
posterior measure of the above inverse problem is well-defined:
\begin{theorem} \label{t:eda} 
Assume that $f \in L^2(0,T,\cH^s)$ with $s>0$. 
Consider the Eulerian data assimilation problem described above. Define a Gaussian measure $\mu_0$ on $\cH$, with mean 
$u_b$ and covariance operator $\beta\,A^{-\alpha}$ for any $\beta>0$ and $\alpha>1.$ If 
$u_b \in \cH^{\alpha},$ then the probability measure $\mu(du)=
\bbP(du|y)$ is absolutely continuous with respect to $\mu_0$
with Radon-Nikodym derivative 
\begin{eqnarray}\label{eq:RNeda} 
\frac{d\mu}{d\mu_0}(u) \propto \exp \left ( -\frac{1}{2\gamma^2}|y - 
\cG(u)|^2_\Sigma \right ).\\
\nonumber 
\end{eqnarray} 
\end{theorem} 

We now define an approximation $\mu^N$ to $\mu$ given by
\eqref{eq:RNeda}. The approximation is made by
employing the Galerkin approximations of $v$ to 
define an approximate $\cG$.
The Galerkin approximation of $v$, $v^N$, is the solution of 
\begin{align}\label{eq:gls}
\frac{\ud v^N}{\ud t}+\nu Av^N+P^N\,B(v^N,v^N)=P^N \psi,\quad v^N(0)=P^N u,
\end{align}
with $P^N$ as defined in the previous section.
Let 
$$
\cG^N(u)=\big(v^N(x_1,t),\dots,v^N(x_K,t)\big)^T
$$
and then consider the approximate prior measure $\mu^N$ defined
via its Radon-Nikodym derivative with respect to $\mu_0$: 
\begin{equation}\label{eq:eumu}
\frac{d \mu^N}{d \mu_0}\propto \exp\left(-\frac{1}{2\gamma^2}|y - 
\cG^N(u)|^2_\Sigma \right ).
\end{equation}
Our aim is to show that $\mu^N$ converges to $\mu$ in the
Hellinger metric. 
Unlike the examples in the previous two sections
we are unable to obtain sufficient control on the
dependence of the error constant on $u$ in the
forward error bound to enable application
of Theorem \ref{t:wp2}; hence we employ Theorem \ref{t:wp2f}.
In the following lemma we obtain a bound on $\|v(t)-v^N(t)\|_{L^\infty(\Omega)}$ 
and therefore on $|\cG(u)-\cG^N(u)|$. Following the
statement of the lemma, we state and prove the
basic approximation theorem for this section.
The proof of the lemma is given 
after the statement and proof
of the approximation theorem for the posterior probability
measure.

\begin{lemma}\label{l:eudiff}
Let $v^N$ be the solution of the Galerkin system (\ref{eq:gls}).
For any $t>t_0$
$$
\|v(t)-v^N(t)\|_{L^\infty(\Omega)}\,\le\, C(\|u\|,t_0)\,\psi(N),
$$
where $\psi(N)\to 0$ as $N\to\infty$.

\end{lemma} 

The above lemma
leads us to the following convergence result for $\mu^N$:

\begin{theorem}
Let $\mu^N$ be defined according to (\ref{eq:eumu}) and let
the assumptions of Theorem \ref{t:eda} hold. Then 
$$
\dhh(\mu,\mu^{N})\to 0
$$
as $N\to\infty$.
\end{theorem}

\begin{proof}
We apply Theorem \ref{t:wp2f} with $X=\cH.$ Assumption \ref{ass:g1}
(and hence Assumption \ref{ass:1}) is established in Lemma 3.1
of \cite{CDRS08}. 
By Lemma \ref{l:eudiff}
$$
|\cG(u)-\cG^N(u)|\le K\psi(N)
$$
with $K=K(\|u\|)$ and $\psi(N)\to 0$ as $N\to 0$. 
Therefore the result follows by Theorem \ref{t:wp2f}.
\end{proof}

\begin{proof}[of Lemma \ref{l:eudiff}]
Let $e_1=v-P^N v$ and $e_2=P^N v-\vN$. 
Applying $P^N$ to (\ref{eq:nsef}) yields
$$
\frac{\ud P^N v}{\ud t}+\nu A P^N v+P^N B(v,v)=P^N \psi.
$$
Therefore $e_2=P^N v-\vN$ satisfies
\begin{equation}\label{eq:e2}
\frac{\ud e_2}{\ud t}+\nu A e_2=P^N B(e_1+e_2,v)+P^N B(\vN,e_1+e_2),\quad e_2(0)=0.
\end{equation}
Since for any and for $m>l$  
\begin{equation}\label{eq:e1decay}
\|e_1\|_{l}^2\le \frac{1}{N^{2(m-l)}}\|v\|_{m}^2,
\end{equation}
we will obtain an upper bound for $\|e_2\|_{1+l}$, $l>0$, in terms of the Sobolev norms of $e_1$ and then use the embedding $\cH^{1+l}\subset L^\infty$ to conclude the result of the lemma.\\

Taking the inner product of (\ref{eq:e2}) with $e_2$, and noting that 
$P^N$ is self-adjoint and $P^N e_2=e_2$ and $(B(v,w),w)=0$, we obtain
\begin{align*}
\frac 1 2\ddt\|e_2\|^2+\nu\|De_2\|^2
&= (B(e_1+e_2,v),e_2)+(B(\vN,e_1),e_2)\\
&\le c\|e_1\|^{1/2}\|e_1\|_1^{1/2}\|v\|_1\|e_2\|^{1/2}\|e_2\|_1^{1/2}
        +c\|e_2\|\,\|v\|_1\,\|e_2\|_1\\
&\quad        +c\|\vN\|^{1/2}\|\vN\|_1^{1/2}\|e_1\|_1\|e_2\|^{1/2}\|e_2\|_1^{1/2}\\
&\le c\|e_1\|^2\,\|e_1\|_1^2+c\|v\|_1^2\,\|e_2\|
       +c\|e_2\|^2\,\|v\|_1^2\\
&\quad  +c\|v^N\|\,\|v^N\|_1\,\|e_1\|_1+c\|e_1\|_1\,\|e_2\|+\frac{\nu}{2}\|e_2\|_1^2
\end{align*}
Therefore
\begin{align*}
\frac{\ud}{\ud t}(1+\|e_2\|^2)
&+\nu\,\|De_2\|^2\le\\ 
&c\,(1+\|v\|_1^2)\,(1+\|e_2\|^2)+c(1+\|e_1\|^2)\,\|e_1\|_1^2
     +c\,\|\vN\|\,\|\vN\|_1\,\|e_1\|_1
\end{align*}
which gives
\begin{align*}
\|e_2(t)\|^2+\nu\int_0^t\|De_2\|^2
&\le c\,\beta(t) \lt( 1+\int_0^t \|\vN\|^2\,\|\vN\|_1^2\,\ud\tau \right)
      \,\int_0^t \|e_1\|_1^2\,\ud\tau\\
&\;+ c\,\beta(t)\,\int_0^t (1+\|e_1\|^2)\,\|e_1\|_1^2\,\ud\tau.
\end{align*}
with 
$$
\beta(t)=\exp\lt(  c\int_0^t 1+\|v\|_1^2\,\ud\tau \right).
$$
Hence
\begin{align}\label{eq:e2-L2}
\|e_2(t)\|^2+\nu\int_0^t\|De_2\|^2 \le c(1+\|u\|^4)\,\mathrm{e}^{c+c\|u\|^2}
\,\int_0^t (1+\|e_1\|^2)\,\|e_1\|_1^2\,\ud\tau.\\
\nonumber
\end{align}

%%%%%%%%%%%%%%%%%%%%%%%%%%%%%%
%\noindent{\bf Estimating $\|e_2(t)\|_s$, for $s<1$:}\\

To estimate $\|e_2(t)\|_s$ for $s<1$, we take the inner product of (\ref{eq:e2}) with $A^{s}e_2$, $0<s<1$ and write
\begin{align*}
\frac{1}{2}\ddt\|e_2\|_{s}^2+\nu\|e_2\|_{1+s}^2 
  \le |\big(((e_1+e_2)\cdot \nabla) v,{A}^{s}e_2\big)|+
      |\big((\vN\cdot \nabla) (e_1+e_2),{A}^{s}e_2\big)|.
\end{align*}
Using 
$$
|\big((u\cdot \nabla) v,{A}^{s}w\big)|\,\le\,c\|u\|_s\,\|v\|_1\,\|w\|_{1+s}
$$
and Young's inequality we obtain
\begin{align*}
\ddt\|e_2\|_{s}^2+\nu\|e_2\|_{1+s}^2 
&\,\le\, c\,(\|e_1\|_s^2+\|e_2\|_s^2)\,\|v\|_1^2
+c\,\|\vN\|_s^2\,(\|e_1\|_1^2+\|e_2\|_{1}^2).
\end{align*}
Now integrating with respect to $t$ over $(t_0,t)$ with $0<t_0<t$ we can write
\begin{align*}
\|e_2(t)\|_{s}^2+\nu\int_{t_0}^t\|e_2\|_{1+s}^2\,\ud\tau 
&\,\le\, \|e_2(t_0)\|_{s}^2+c\,\sup_{\tau\ge t_0}\|v(\tau)\|_1^2 
            \int_0^t \|e_1\|_s^2+\|e_2\|_s^2\ud\tau\\
&\quad +c\,\sup_{\tau\ge t_0}\|\vN(\tau)\|_s^2 
            \int_0^t \|e_1\|_1^2+\|e_2\|_1^2\,\ud\tau.
\end{align*}
Therefore since for $s\le 1$ and $t\ge t_0$
$$
\|v(t)\|_s^2 \le \frac{c(1+\|u\|^2)}{t_0^s},
$$
and noting that the same kind of decay bounds that hold for $v$ can be shown similarly for $\vN$ as well, we have
$$
\|e_2(t)\|_{s}^2+\nu\int_{t_0}^t\|e_2\|_{1+s}^2\,\ud\tau 
\,\le\, \|e_2(t_0)\|_{s}^2 
+ \frac{c}{t_0}(1+\|u\|^6)\mathrm{e}^{c+c\|u\|^2}\,\int_0^t (1+\|e_1\|^2)\,\|e_1\|_1^2\,\ud\tau.
$$
Integrating the above inequality with respect to $t_0$ in $(0,t)$ 
we obtain
\begin{align}\label{eq:decay<1}
\|e_2(t)\|_{s}^2+\nu\int_{t_0}^t\|e_2\|_{1+s}^2\,\ud\tau 
\le \frac{c}{t_0}(t_0+1+\|u\|^6)\,\int_0^t (1+\|e_1\|^2)\,\|e_1\|_1^2\,\ud\tau
\end{align}
for $t>t_0$.\\

%%%%%%%%%%%%%%%%%%%%%%%%%%%%%%%
%\noindent{\bf Estimating $\|e_2(t)\|_s$, for $s>1$:}\\

Now we estimate $\|e_2(t)\|_s$ for $s>1$. Taking the inner product of (\ref{eq:e2}) with $A^{1+l}e_2$, $0<l<1$, we obtain
\begin{align*}
\frac{1}{2}\ddt\|e_2\|_{1+l}^2+\nu\|e_2\|_{2+l}^2 
&  \le |\big(((e_1+e_2)\cdot \nabla) v,{A}^{1+l}e_2\big)|\\
& +      |\big((\vN\cdot \nabla) (e_1+e_2),{A}^{1+l}e_2\big)|.
\end{align*}
Since (see \cite{CDRS08})
\begin{align*} 
\left( (u\cdot\nabla)v, A^{1+l}w \right) 
\le 
c\,\|u\|_{1+l}\,\|v\|_1\,\|w\|_{2+l}+c\,\|u\|_l\,\|v\|_2\,\|w\|_{2+l}
\end{align*} 
and using Young's inequality, we can write
\begin{align*}
\ddt\|e_2\|_{1+l}^2+\nu\|e_2\|_{2+l}^2 
&\le c\,\|e_1\|_{1+l}^2\|v\|_1^2 + c\, \|e_1\|_{l}^2\|v\|_2^2\\
&\quad + c\,\|e_2\|_{1+l}^2\|v\|_1^2 + c\, \|e_2\|_{l}^2\|v\|_2^2\\
&\quad + c\,\|\vN\|_{1+l}^2\|e_1\|_1^2 + c\, \|\vN\|_{l}^2\|e_1\|_2^2\\
&\quad + c\,\|\vN\|_{1+l}^2\|e_2\|_1^2 + c\, \|\vN\|_{l}^{2/l}\|e_2\|_{1+l}^2.
\end{align*}
Now we integrate the above inequality with 
respect to $t$ and over $(t_0/2+\sigma,t)$ with $0<t_0<t$ and $0<\sigma<t-t_0/2$
and obtain (noting that $\|\vN\|_s\le \|v\|_s$ for any $s>0$)
\begin{align*}
\|e_2(t)\|_{1+l}^2\,
&\le \|e_2(t_0/2+\sigma)\|_{1+l}^2
         +\sup_{\tau\ge t_0/2}\|v(\tau)\|_1^2 
                   \int_{t_0/2+\sigma}^t\|e_1\|_{1+l}^2+\|e_2\|_{1+l}^2\,\ud\tau\\
&\quad + \sup_{\tau\ge t_0/2}(\|e_1(\tau)\|_l^2+ \|e_2(\tau)\|_l^2)
                   \int_{t_0/2+\sigma}^t\|v\|_2^2\,\ud\tau\\
&\quad + \sup_{\tau\ge t_0/2}(\|e_1(\tau)\|_1^2+ \|e_2(\tau)\|_1^2)
                   \int_{t_0/2+\sigma}^t\|\vN\|_{1+l}^2\,\ud\tau\\
&\quad +\sup_{\tau\ge t_0/2}(1+\|\vN(\tau)\|_l^{2/l} )
                   \int_{t_0/2+\sigma}^t\|e_1\|_{2}^2+\|e_2\|_{1+l}^2\,\ud\tau.
\end{align*}
We have, for $s>1$ and $t>t_0$, (\cite{CDRS08})
$$
\|v(t)\|_s^2\le \frac{c(1+\|u\|^4)}{t_0^s}.
$$
Therefore using (\ref{eq:decay<1}) and (\ref{eq:e1decay}) we conclude that
\begin{align*}
\|e_2(t)\|_{1+l}^2
&\le\, \|e_2(t_0/2+\sigma)\|_{1+l}^2 \\
& +C_p(\|u\|)\lt(\frac {1} { N^{2(m-l)} \, t_0^{1+m} }
         +\frac {1} { t_0^{1+l} }\,\int_0^t (1+\|e_1\|^2)\,\|e_1\|_1^2\,\ud\tau
            +\frac {1} { N^{2(r-1)} \, t_0^{1+r} }  \right)
\end{align*}
with $r>1$ and where $C_p(\|u\|)$ is a constant depending on polynomials of $\|u\|$.
Integrating the above inequality with respect to $\sigma$ over $(0,t-t_0/2)$ we obtain
\begin{align*}
\|e_2(t)\|_{1+l}^2 
&\le 
  C_p(\|u\|) \lt( \frac{1}{t_0^{1+l} }+\frac{1}{t_0^{2+l}} \right)\,\int_0^t (1+\|e_1\|^2)\,\|e_1\|_1^2\,\ud\tau\\
&\quad    +C_p(\|u\|)\lt(\frac {1} { N^{2(m-l)} \, t_0^{2+m} }
            +\frac {1} { N^{2(r-1)} \, t_0^{2+r} }  \right).
\end{align*}
Now to show that $\|e_1\|^2+\int_0^t \|e_1\|_1^2\,\ud\tau\to 0$ as $N\to\infty$, we note that $e_1$ satisfies
\begin{align*}
\frac{1}{2}\ddt\|e_1\|^2+\nu\|De_1\|
&\le\|(I-\PN)f\|\,\|e_1\|+\|(B(v,v),e_1)\|\\
&\le\|(I-\PN)f\|\,\|e_1\|+\|v\|^{1/2}\,\|Dv\|^{3/2}\,\|e_1\|^{1/2}\,\|De_1\|^{1/2}\\
&\le\|(I-\PN)f\|\,\|e_1\|+c\,\|v\|^{2/3}\,\|Dv\|^{2}\,\|e_1\|^2+\frac{\nu}{2}\|De_1\|^2.
\end{align*}
Therefore
$$
\ddt\|e_1\|^2+\nu\|De_1\|\le\,c\,\|(I-\PN)f\|^2+c\,(1+\|v\|^{2/3}\,\|Dv\|^{2})\,\|e_1\|^2
$$
and after integrating, we get
$$
\|e_1\|^2+\int_0^T\|e_1\|_1\,\ud\tau\,\le \exp(1+C_p(\|u\|))\,\lt(\|e_1(0)\|^2
+\int_0^T\|(I-\PN)f\|^2 \right)\,\ud\tau.
$$
Since $f\in L^2(0,T;\cH)$, the above integral tend to zero as $N\to\infty$ and the result follows.
\end{proof}

\section{Conclusions}

In this paper we have studied the approximation of
inverse problems which have been regularized by means
of a Bayesian formulation. We have developed a general
approximation theory which allows for the transfer
of approximation results for the forward problem
into approximation results for the inverse problem.
The theory clearly separates analysis of the forward
problem, in which no probabilistic methods are
required, and the probabilistic framework for the
inverse problem itself:
it is simply necessary that the
requisite bounds and approximation properties
for the forward problem hold in a space
with full measure under the prior.
Indeed the approximation theory may be
seen to place constraints on the prior,
in order to ensure the desired robustness.

In applications there are two sources of  
error when calculating expectations
of functions of infinite dimensional random
variables: the error which we provide an analysis
for in this paper, namely the approximation of
the measure itself in a finite dimensional
subspace, together with the error
incurred through calculation of expectations.
The latter can be undertaken by 
Markov chain-Monte Carlo (MCMC) methods, or quasi Monte Carlo methods.
The two sources of error must be balanced in order
to optimize computational cost. 

We have studied three specific applications, all
concerned with determining the initial condition 
of a dissipative PDE, from observations of various
kinds, at positive times. However the general
approach is applicable to a range of inverse
problems for functions when formulated in a Bayesian fashion.
The article \cite{Stuart10} overviews many applications
from this point of view. Furthermore we have
limited our approximation of the underlying
forward problem to spectral methods. However we
anticipate that the general approach will be useful
for the analysis of other spatial approximations based on
finite element methods, for example, and to approximation
errors resulting from time-discretization; indeed it
would be interesting to carry out analyses for such
approximations.

It is important to realize that
new approaches to the computation of expectations
against measures on infinite dimensional spaces
are currently an active area of research
in the engineering community \cite{SpG89,GSp03} and that
a numerical analysis of this area is being systematically
developed \cite{ST06,TS07}. That work is primarily
concerned with approximating measures which are the
push forward, under a nonlinear map, of a simple
measure with product strcuture, such as a Gaussian measure;
in contrast the inverse problem setting which we
study here is concerned with the approximation of
non-Gaussian measures whose Radon-Nikodym derivative
is defined through a related nonlinear map.
It would be interesting
to combine the approaches
in \cite{SpG89,ST06,TS07} and related literature with
the approximation theories described in this paper. 
For example that work could be used to develop
cheap approximations to the forward map $\cG$ 
thereby accelerating
MCMC-based sampling methods.

\vspace{0.1in}

\noindent{\bf Acknoweldgements} The authors are grateful to
the EPSRC, ERC and ONR for financial support.

\vspace{0.1in}

%\newpage

\bibliography{mybib}
\bibliographystyle{plain}

%\newpage

\begin{appendix}
\section{Analytic Semigroups and Probability}

We collect together some basic facts concerning analytic
semigroups and probability
required in the main body of
the article.
First we state the well-known Gronwall inequality
in the form in which we will use it\footnote{
See {\tt http://en.wikipedia.org/wiki/Gronwall's$\textunderscore$inequality}}

\begin{lemma}\index{Gronwall inequality}\label{lem:gron}
Let $I=[c,d)$ with $d \in (c,\infty].$ Assume that
$\alpha,u \in C(I;\bbR^+)$ and that there is $\lambda<\infty$
such that, for all intervals $J\subseteq I$,
$\int_J \beta(s)ds<\lambda.$ 
If
$$u(t) \le \alpha(t)+\int_c^t \beta(s)u(s)ds,\quad t\in I,$$
then
$$u(t) \le \alpha(t)+\int_c^t \alpha(s)\beta(s)
\exp\Bigl(\int_s^t \beta(r)dr\Bigr)ds,\quad t\in I.$$
In particular, if $\alpha(t)=u+2at$ is positive in $I$
and $\beta(t)=2b$ then
$$u(t) \le \exp(2bt)u+\frac{a}{b}\Bigl(\exp(2bt)-1 \Bigr),
\quad t\in I.$$
Finally, if $c=0,$ and $0<\alpha(t) \le K$ in $I,$ then
$$u(t) \le K+K\lambda\exp(\lambda), \quad t\in I.$$
\end{lemma}

Throughout this article $A$ denotes either
the Laplacian on a smooth, bounded domain in $\bbR^d$
with Dirichlet boundary conditions (section \ref{sec:heat})
or the Stokes operator on $\bbT^2$ (sections \ref{sec:lag}
and \ref{sec:eul}). In both cases $A$ is
a self-adjoint  positive operator $A$, 
densely defined on a Hilbert space $\cH$, and the generator
of an analytic semigroup.
We denote by $\{(\phi_k,\lambda_k)\}_{k\in \bbK}$ a complete
orthonormal set of eigenfunctions/eigenvalues for $A$ in $\cH.$
We then define fractional powers of $A$ by 
\begin{equation}
\label{eq:carf2}
A^{\alpha}u=\sum_{k\in\bbK} \lambda_k^{\alpha} \langle u,\phi_k\rangle \phi_k.
\end{equation}
For any $s\in \bbR$ we define the Hilbert 
spaces\index{Hilbert space} $\cHs$ by 
\begin{equation}
\label{eq:Hs}
\cHs=\{u: \sum_{k \in \bbK} \lambda_{k}^s |\langle u,\phi_k\rangle|^2<\infty\}.
\end{equation}
The norm in $\cHs$ is denoted by $\|\cdot\|_{s}$ and is given by 
$$\|u\|_s^2=\sum_{k \in \bbK} \lambda_{k}^s |\langle u,\phi_k\rangle|^2.$$
Of course $\cH^0=\cH$.
If $s > 0$ then these spaces are contained in $\cH$, but for
$s<0$ they are larger than $\cH$.
It follows that the domain of $A^{\alpha}$ is 
$\cH^{2\alpha};$
the image of $A^{-\alpha}$ is 
$\cH^{2\alpha}.$

Now consider the Hilbert-space
valued ODE
\begin{equation}
\frac{dv}{dt}+A v=f, \quad v(0)=u. \label{eq:par}
\end{equation}
We state some basic results in this area, provable
by use of the techniques in \cite{pagy}, for example,
or by direct calculation using the eigenbasis for $A.$
For $f=0$ the solution 
$v \in C([0,\infty),\cH)\cap C^1((0,\infty),D(A))$ and
\begin{equation}
\|v\|_{s}^2 \le Ct^{-(s-l)} \|u\|_{l}^2, \quad \forall t \in (0,T].
\label{t:smooth}
\end{equation}
If $f \in C([0,T],\cH^\gamma)$ for some
$\gamma \ge 0$, then \eqref{eq:par} 
has a unique mild solution
$u \in C([0,T];\cH)$ and, for $0 \le \ell <\gamma+2$,
\begin{equation}
\label{t:thirteen}
\|v(t)\|_s \le C\Bigl(\frac{\|u\|_l}{t^{(s-l)/2}}+
\|f\|_{C([0,T],\cH^\gamma)}\Bigr)
\end{equation}
for $s\in [\ell,2+\gamma).$

It central to this paper to estimate
the distance between two probability measures. 
To this end we introduce two useful metrics 
on measures: the {\em total variation distance}\index{total
variation distance} and
the {\em Hellinger distance}\index{Hellinger distance}.
We discuss the relationships between the metrics
and indicate how they may be used to estimate differences
between expectations of random variables under two different
measures.

Assume that we have two probability measures $\mu$ and $\mu'$,
both absolutely continuous with respect to the
same reference measure $\nu$. The following
defines two concepts of distance between $\mu$ and $\mu'$.

\begin{defn} \label{def:tvd}
The {\em total variation distance}\index{total variation
distance} between
$\mu$ and $\mu'$ is
$$\dtv(\mu,\mu')=\frac12 \int \Bigl|
\frac{d\mu}{d\nu}-\frac{d\mu'}{d\nu}\Bigr|
d\nu.$$
The {\em Hellinger distance}\index{Hellinger distance} between
$\mu$ and $\mu'$ is
$$\dhh(\mu,\mu')=\sqrt{\Bigl(\frac12 \int \Bigl(
\sqrt\frac{d\mu}{d\nu}-\sqrt\frac{d\mu'}{d\nu}\Bigr)^2
d\nu\Bigr)}.$$
\end{defn}

Both distances are invariant under the choice of
$\nu$ in that they are unchanged if a different reference
measure, with respect to which $\mu$ and $\mu'$ are
absolutely continuous, is used. Furthermore, it
follows from the definitions that
$\dtv(\mu,\mu')\in (0,1)$ and $\dhh(\mu,\mu')\in (0,1).$
The Hellinger and total variation distances are related
as follows\cite{GS02}\footnote{Note that different
normalization constants are sometimes used in the definitions
of distance.}:
\begin{equation}
\label{eq:tvhrel}
\frac{1}{\sqrt{2}}\dtv(\mu,\mu') \le \dhh(\mu,\mu') \le
\dtv(\mu,\mu')^{\frac12}.
\end{equation}

The Hellinger distance is particularly useful
for estimating the difference between expectation
values of functions of random variables under
different measures. This is illustrated in the
following lemma:

\begin{lemma} \label{l:tvhell2}
Assume that two measures $\mu$ and $\mu'$ on a Banach
space $\Bigl(X,\|\cdot\|_{X}\Bigr)$
are both absolutely continuous with respect to a measure
$\nu.$ Assume also that $f:X \to Z$, where 
$\Bigl(Z,\|\cdot\|\Bigr)$ is a Banach space, 
has second moments with respect to
both $\mu$ and $\mu'$. Then
$$\|\bbE^{\mu}f-\bbE^{\mu'}f\| \le 2\Bigl(\bbE^{\mu}\|f\|^2+
\bbE^{\mu'}\|f\|^2\Bigr)^{\frac12}\dhh(\mu,\mu').$$
Furthermore, if $\Bigl(Z,\langle \cdot,\cdot \rangle\Bigr)$ is
a Hilbert space and $f:X \to Z$ has fourth moments then
$$\|\bbE^{\mu}f\otimes f-\bbE^{\mu'}f\otimes f\| \le 2\Bigl(\bbE^{\mu}\|f\|^4+
\bbE^{\mu'}\|f\|^4\Bigr)^{\frac12}\dhh(\mu,\mu').$$
\end{lemma}

\begin{proof} We have
\begin{align*}
\|\bbE^{\mu}f-\bbE^{\mu'}f\|
&\le \int \|f\| \Bigl|\frac{d\mu}{d\nu}-\frac{d\mu'}{d\nu}\Bigr|d\nu\\
& \le  \int \Bigl(\frac{1}{\sqrt 2}\Bigl|\sqrt{\frac{d\mu}{d\nu}}-
\sqrt{\frac{d\mu'}{d\nu}}\Bigr|\Bigr)\Bigl(
\sqrt{2}\|f\|\Bigl|\sqrt{\frac{d\mu}{d\nu}}+
\sqrt{\frac{d\mu'}{d\nu}}\Bigr|\Bigr) d\nu\\
& \le \sqrt{\Bigl(\frac12 \int \Bigl(
\sqrt\frac{d\mu}{d\nu}-\sqrt\frac{d\mu'}{d\nu}\Bigr)^2
d\nu\Bigr)}
\sqrt{\Bigl(2 \int \|f\|^2\Bigl(
\sqrt\frac{d\mu}{d\nu}+\sqrt\frac{d\mu'}{d\nu}\Bigr)^2
d\nu\Bigr)}\\
& \le \sqrt{\Bigl(\frac12 \int \Bigl(
\sqrt\frac{d\mu}{d\nu}-\sqrt\frac{d\mu'}{d\nu}\Bigr)^2
d\nu\Bigr)}
\sqrt{\Bigl(4 \int \|f\|^2\Bigl( 
\frac{d\mu}{d\nu}+\frac{d\mu'}{d\nu}\Bigr)
d\nu\Bigr)}\\
&= 2\Bigl(\bbE^{\mu}\|f\|^2+
\bbE^{\mu'}\|f\|^2\Bigr)^{\frac12}\dhh(\mu,\mu')
\end{align*}
as required. 

The proof for $f\otimes f$ follows from the following
inequalities, and then arguing similarly to the case
for the norm of $f$:
\begin{align*}
\|\bbE^{\mu}f\otimes f-\bbE^{\mu'}f\otimes f\|
&=\sup_{\|h\|=1}
\|\bbE^{\mu}\langle f,h\rangle f-\bbE^{\mu'}\langle f,h\rangle f\|\\
&\le \int \|f\|^2 \Bigl|\frac{d\mu}{d\nu}-\frac{d\mu'}{d\nu}\Bigr|d\nu.
\end{align*}
\end{proof}

Note, in particular, that choosing $X=Z$, and with $f$ chosen
to be the identity mapping, we deduce that the differences
in mean and covariance operators under two measures are bounded
above by the Hellinger distance between the two measures.

The following Fernique Theorem (see \cite{DapZab92}, Theorem 2.6) 
will be used repeatedly:

\begin{theorem} \label{thm:fern} Let $x \sim \mu=\cN(0,\cC)$
where $\mu$ is a Gaussian measure on Hilbert space $H$.
Assume that $\mu_0(X)=1$ for some Banach space $\Bigl(X,
\|\cdot\nx\Bigr)$ with $X \subseteq H.$ Then
there exists $\alpha>0$ such that
$$\int_{\cH} \exp\bigl(\alpha \|x\nx^2\bigr)\mu(dx)<\infty.$$
\end{theorem}

The following regularity properties of Gaussian random
fields will be useful to us; the results may be proved
by use of the Kolmogorov continuity criterion, together
with the Karhunen-Loeve expansion (see \cite{DapZab92},
section 3.2):

\begin{lemma} \label{lem:greg2}
Consider a Gaussian measure $\mu=\cN(0,\igmas)$
with $\igmas=\beta A^{-\alpha}$ where $A$ is as defined
earlier in this Appendix A. 
Then $u\sim\mu$ is almost surely
$s-$H\"{o}lder continuous for any exponent $s<
\min\{1,\alpha-\frac{d}{2}\}$
and $u \in \cH^s$, $\mu-$almost surely,
for any $s<\alpha-\frac{d}{2}.$
\end{lemma}

\end{appendix}

\end{document}